\documentclass[11pt,a4paper]{article}
\usepackage{customstyle}
\usepackage{authblk}
\usepackage{color,soul}

\usepackage{lineno}

\title{\textbf{SINDybrid: automatic generation of hybrid models for dynamic systems}}
\author[1]{Ulderico Di Caprio\thanks{Corresponding author: \texttt{ulderico.dicaprio@kuleuven.be}}}
\author[1]{M. Enis Leblebici\thanks{Corresponding author: \texttt{muminenis.leblebici@kuleuven.be}}}
\affil[1]{Center for Industrial Process Technology, Department of Chemical Engineering, KU Leuven, Agoralaan Building B, 3590 Diepenbeek, Belgium}
\date{\vspace{-5ex}}

\begin{document}
\maketitle
\vspace{-1.5em} 

\begin{abstract}
Hybrid modelling enhances the accuracy and predictive capability of dynamic models by integrating first principles with data-driven methods, effectively mitigating epistemic uncertainties inherent in mechanistic approaches. However, hybrid model construction remains complex, typically requiring expert knowledge to identify model epistemic uncertainty and select suitable machine-learning components to capture it. This complexity limits broader adoption in research and industry. We introduce SINDybrid, an automated algorithm designed to streamline hybrid model development for dynamic systems. SINDybrid employs a mixed-integer linear programming (MILP) approach to systematically identify epistemic uncertainty sources and compensate them using optimally selected data-driven components. For broader accessibility and reproducibility, we provide SINDybrid as an open-source Python library. SINDybrid was validated through three case studies: a catalytic reaction system, a continuous fermentation process, and a Lotka-Volterra oscillator, each showcasing different epistemic uncertainties. The robustness of the algorithm is tested against varying experimental conditions, including measurement noise (up to 20\%), limited training batches (minimum 2), and sparse temporal sampling (minimum 5 samples per experiment). Results demonstrate a consistent ability of the SINDybrid approach to produce accurate hybrid models, achieving R² scores above 0.85 on validation data. Additionally, the algorithm effectively identifies uncertainty locations within system dynamics under challenging scenarios. This work highlights SINDybrid potential as a versatile, automated solution for hybrid modelling, significantly reducing the barriers to adopting hybrid methods in complex dynamic systems across scientific research and industrial applications.
\end{abstract}

\keywords{Hybrid modelling; SINDy; Data-driven models; Epistemic uncertainty; MILP}

\section{Introduction}
Digital tools are transforming the chemical industry by changing how we control and optimise processes to achieve greater efficiency and sustainability. Modelling and simulation are utilised at various stages of chemical production, from R\&D to manufacturing and operations \citep{Sansana_Joswiak_Castillo_Wang_Rendall_Chiang_Reis_2021, Chiang_etal_2022}. Traditionally, process models are developed using first-principles modelling (FPM), which involves formulating equations that describe system behaviour based on a mechanistic understanding. FPMs provide advantages such as physical interpretability and insights across multiple scales. However, creating accurate FPMs requires extensive knowledge of all relevant phenomena to minimise epistemic uncertainty of the description and produce accurate and reliable predictions. Addressing this uncertainty often necessitates comprehensive, system-specific experimentation, which can be both costly and time-consuming.
In contrast, data-driven modelling (DDM) and machine learning (ML) approaches do not require detailed mechanistic understanding. These models rely on statistical representations with optimisable parameters, which are adjusted or tuned using experimental data and use experimental data to map input-output relationships. Because of this, DDM and ML are increasingly adopted in chemical engineering tasks \cite{Schweidtmann_etal_2021_review}. Nonetheless, DDMs inherently lack adherence to physical laws (e.g., conservation principles and transport equations), which can compromise robustness and extrapolation capabilities. Moreover, they typically require large, high-quality datasets and offer limited means to incorporate existing domain knowledge.
Hybrid modelling (HM) combines the strengths of FPM and DDM by embedding data-driven components within mechanistic frameworks \citep{vonStosch_etal_2014, Schweidtmann_etal_2024_review}. In these models, DDMs are typically used to compensate for epistemic uncertainty in mechanistic equations, enabling the capture of complex, poorly understood behaviours while retaining consistency with physical laws. Despite their potential, HMs remain difficult to construct. Their development demands interdisciplinary expertise, spanning process modelling, ML, and domain-specific knowledge, making them challenging to design and validate. In this scenario, an automated framework that automates HM construction and facilitates its broader adoption in the chemical industry is beneficial. Ideally, such a framework would allow users to provide only mechanistic knowledge while identifying the equation running under epistemic uncertainty and compensating for it, selecting and training a DDM.
Chemical processes are typically governed by multiple interacting physical phenomena (e.g., heat and mass transfer, reaction kinetics, fluid dynamics) and often exhibit temporal dependencies. As a result, FPMs often comprise systems of ordinary differential equations (ODEs) or differential-algebraic equations (DAEs), each describing a specific physical mechanism. These equations contain parameters (e.g., reaction rate constants) that must be estimated from literature or experiments. Epistemic uncertainty arises when certain mechanisms are not fully comprehended or captured in the mechanistic model. Developing an effective HM, therefore, requires identifying the model regions (i.e., equations or parameters) most affected by uncertainty and selecting appropriate data-driven components to correct them. Automating this process involves two key tasks: (1) identifying sources of epistemic uncertainty and (2) selecting suitable DDM structures for compensation. While much of the existing literature focuses on the second task \citep{Willis_etal_2017, Narayanan_etal_2022, Wilson_Sahinidis_2017, DeCarvalho_etal_2024, Zhang_Savage_Cho_2020}, the first remains underexplored.
Sparse identification is a modern approach that enables the discovery of governing equations from data in dynamic systems, offering a compact and interpretable alternative to traditional DDMs. These approaches use a predefined library of candidate functions and apply L1-regularised regression (i.e., LASSO) to identify the most relevant terms from data. The L1 regularisation promotes sparsity by driving irrelevant terms to zero, yielding compact and interpretable models. The Sparse Identification of Nonlinear Dynamical Systems (SINDy) framework, introduced by Brunton et al.\cite{SINDy_paper, SINDyC}, formalised this idea for both autonomous and forced systems. Subsequent developments have extended SINDy to handle nonlinear parameters \cite{Viknesh_Tatari_Arzani_2025}, nested symbolic functions \cite{SINDy_ANN_2025}, hybrid dynamics \cite{Mangan_Askham_Brunton_Kutz_Proctor_2019}, tensor decomposition for scalability \cite{Gelss_etal_2019} and Bayesian inference under uncertainty \cite{Hirsh_etal_2021}. Due to its versatility and accuracy, SINDy approaches have been applied to various chemical engineering tasks, including system modelling \cite{Zhang_Savage_Cho_2020, Abdullah_Alhajeri_Christofides_2022, Abdullah_Christofides_2023a, Abdullah_Christofides_2023b} and reaction mechanism discovery \cite{Hoffmann_Fröhner_Noé_2019}.
However, most existing SINDy-based methods ignore prior mechanistic knowledge. They attempt to rediscover the entire system of equations from data, regardless of whether some components are already known. This leads to inefficient use of data and limits the interpretability of the resulting model. A methodology that selectively applies SINDy only to uncertain parts of the model, while preserving validated mechanistic equations, would significantly improve both data efficiency and hybrid model reliability.
This paper addresses this gap by introducing SINDybrid, a novel SINDy-based framework for automated hybrid model generation. SINDybrid selectively integrates data-driven corrections into mechanistic models. It automates the identification of uncertain components, produces symbolic and interpretable corrections to enhance model fidelity, and demonstrates robustness across three benchmark dynamic systems under noisy and sparse conditions. SINDybrid integrates symbolic differential equations representing mechanistic knowledge with sparse identification applied selectively to uncertain components. The method automatically detects epistemic uncertainty and compensates for it using interpretable, symbolic expressions learned from data. This results in hybrid models that are physically consistent, data-efficient, and easy to interpret, validate, and share. To the best of our knowledge, this is the first framework that systematically combines mechanistic modelling and sparse identification to automate hybrid model development in dynamic chemical systems.

\section{Mathematical formulation}
\subsection{Case description and optimisation formulation}
Say a dynamic system with N\textsubscript{S} states whose real behaviour is described by the ODE model
\begin{equation} \label{eq:1}
\frac{d{\boldsymbol{x}}}{dt} = \boldsymbol{f}({\boldsymbol{x}}, \boldsymbol{r}) \space + \space \boldsymbol{h}({\boldsymbol{x}}, \boldsymbol{r})
\end{equation}
where $\boldsymbol{x}$ is a column vector $(N_S,1)$ containing the system states, $\boldsymbol{f}: \mathbb{R}^{N_S} \rightarrow \mathbb{R}^{N_S}$ is the equation set of the FPM containing the available mechanistic information about the system, and $\boldsymbol{h}: \mathbb{R}^{N_S} \rightarrow \mathbb{R}^{N_S}$ is the epistemic uncertainty of the system whose functional form is not available. The function $\boldsymbol{h}$ is not known by the designer and contains the quantitative description of the phenomena occurring in the system, but that are not taken into account in the formulation of $\boldsymbol{f}$ (i.e., the epistemic uncertainty). Therefore, the non-null terms of $\boldsymbol{h}$ are the ones corresponding to the mechanistic model equations running under epistemic uncertainty. $\boldsymbol{r}$ are the further conditions used to characterise the system, but that are constant during an experimental run characterising the experiment(e.g., temperature, catalyst concentration). In the proposed formulation, constructing an HM translates to identifying the functional forms of the terms contained in $\boldsymbol{h}$. Say an experimentation campaign has been executed on the system, and $N_t$ experimental samples were obtained for each experiment. In this condition, each experiment is represented by a matrix $\boldsymbol{x}_{\text{exp}}$ with dimension $(N_t, N_S)$, where each row contains a time-sample and the columns contain the system states. From the matrix $\boldsymbol{x}_{\text{exp}}$, it is possible to numerically calculate the time derivative values for all the $N_S$ system states (i.e., $d{\boldsymbol{x}_{exp}}/dt$). Therefore, employing \eqref{eq:1} one can calculate the numerical values of the $\boldsymbol{h}$ matrix and numerically characterise the model epistemic uncertainty
\begin{equation}\label{eq:2}
{\boldsymbol{h}}_{exp} = \boldsymbol{h}(\boldsymbol{x}_{exp}, \boldsymbol{r}_{exp}) = \frac{d{\boldsymbol{x}_{exp}}}{dt} - \boldsymbol{f}(\boldsymbol{x}_{exp}, \boldsymbol{r}_{exp}).
\end{equation}

In \eqref{eq:2}, $\boldsymbol{h}_{exp}$ is a matrix with known coefficients that can be fitted with the experimental conditions and states using data-driven methodologies, such as SINDy approaches \cite{SINDy_paper}. However, only some of the columns of $\boldsymbol{h}_{exp}$ are non-zero and correspond to the state equations affected by epistemic uncertainty.
Theoretically, classical SINDy can solve this problem thanks to the possibility of selecting relevant features for the problem. However, in practice, numerical differentiation introduces noise, resulting in $\boldsymbol{h}_{exp}$ matrices with non-zero values across all entries. An example of this effect can be found in Figure \ref{fig:7} of this manuscript. 
A potential solution is to increase the regularisation parameter, but this could suppress relevant system dynamics to simplify the model structure, thereby reducing model accuracy. An alternative approach involves constraining the optimisation process while enabling the algorithm to identify the specific equations where deviations occur and identify only the parameters related to them. In this work, we adopt the latter strategy, enforcing constraints that allow the optimiser to determine which equations should include the data-driven correction.
In the SINDy methodology, the matrix $\boldsymbol{h}_{exp}$ is approximated using a library of known functions evaluated on the experimental state values $\boldsymbol{x}_{exp}$. Specifically, SINDy constructs a library matrix $\boldsymbol{X}_L \in \mathbb{R}^{N_t \times N_L}$, where $N_L$ is the number of candidate functions (i.e., the number of functions contained in the library) and $N_t$ is the number of time samples. Each column of $\boldsymbol{X}_L$ includes values of a different candidate function evaluated at the state trajectories.
For instance, the library may include monomials of the system states (e.g., $x_i$, $x_i^2$, $x_i x_j$), trigonometric terms (e.g., $\sin(x_i)$, $\cos(x_j)$), or custom made functions (e.g., $(0.1 + x_1)^{-1}$). Each term is associated with a tunable coefficient $\xi_{ij}$, representing the contribution of the $i$-th function to the model of the $j$-th state. Further information about the library generation and the data-driven hypermodel construction can be found in the SINDy paper \cite{SINDy_paper}.
Though this formulation, the matrix $\boldsymbol{h}_{pred}$ can be expressed as
\begin{equation}\label{eq:3}
\boldsymbol{h}_{pred} = \boldsymbol{X}_L \cdot \boldsymbol{\Xi},
\end{equation}
where $\boldsymbol{h}_{pred}$ is the matrix containing the predicted deviation values, $\boldsymbol{X}_L$ is a matrix containing the values of the functions contained in the library, and $\boldsymbol{\Xi}$, with dimensions $(N_L, N_S)$, is the matrix containing the parameters $\xi$ that need to be identified through the model training. Therefore, with the approach just illustrated, the design of the HM has been converted to an identification problem of the parameters contained in $\boldsymbol{\Xi}$ under the constraint that only a fixed number of columns should have non-null elements. The identification is driven by minimising the loss function containing the sum of the absolute errors between the model prediction and the actual system deviations
\begin{equation}\label{eq:4}
L(\boldsymbol{\Xi}) = \sum_{i,j=0}^{N_t,N_S} \left| \boldsymbol{h}_{pred} - \boldsymbol{h}_{exp} \right|
\end{equation}
and an L1 regularisation can be applied to let the optimiser select only relevant library functions for the analysed case, reshaping the loss function as
\begin{equation}\label{eq:5}
L(\boldsymbol{\Xi}) = \sum_{i,j=0}^{N_t,N_S} \left| \boldsymbol{h}_{pred} - \boldsymbol{h}_{exp} \right| + \lambda_{1,\xi} \sum_{i,j=0}^{N_L,N_S} \left| \xi_{ij} \right|
\end{equation}
where $\lambda_{1,\xi}$ is a regularisation parameter chosen prior the identification execution. Additionally, a further regularisation to the loss function is needed to minimise the amount of active data-driven models (i.e., active terms of $\boldsymbol{h}$), including $\lambda_{1,\delta} \sum_{j=0}^{N_S} \delta_j$ as a further term of the loss function \eqref{eq:5}, resulting into

\begin{equation}\label{eq:5bis}
L(\boldsymbol{\Xi}) = \sum_{i,j=0}^{N_t,N_S} \left| \boldsymbol{h}_{pred} - \boldsymbol{h}_{exp} \right| + \lambda_{1,\xi} \sum_{i,j=0}^{N_L,N_S} \left| \xi_{ij} \right| + \lambda_{1,\delta} \sum_{j=0}^{N_S} \delta_j.
\end{equation}

The value of the parameters $\lambda_{1,\xi}$ and $\lambda_{1,\delta}$ were correlated in order to homogenise the regularisations magnitudes and remove the risk that the optimisation employed many DDM with low values or including only one model with big parameters value. Further information about the derivation of the new regularisation term is reported in the Supplementary Information (SI) section of this manuscript. The final loss function after the regularisation is
\begin{equation}\label{eq:6}
L(\boldsymbol{\Xi}) = \sum_{i,j=0}^{N_t,N_S} \left| \boldsymbol{h}_{pred} - \boldsymbol{h}_{exp} \right| + \lambda_{1,\xi} \left( \sum_{i,j=0}^{N_L,N_S} \left| \xi_{ij} \right| + \max(|ub|, |lb|) \cdot N_L \cdot \sum_{j=0}^{N_S} \delta_j \right)
\end{equation}
Given the loss function, the problem can be solved given the maximum cardinality of the system $K_\alpha$ (i.e., selecting how many columns of $\boldsymbol{\Xi}$ should be non-null a-priori) and the maximum allowed number of active DDMs $K_\delta$. Overall, the optimisation formulation results in

\begin{align}\label{eq:7}
\min_{\boldsymbol{\Xi}, \boldsymbol{A}, \boldsymbol{\Delta}} \quad 
& \sum_{i,j=0}^{N_t,N_S} \left| \boldsymbol{h}_{pred} - \boldsymbol{h}_{exp} \right| 
+ \lambda_{1,\xi} \left( \sum_{i,j=0}^{N_L,N_S} \left| \xi_{ij} \right| 
+ \max(|ub|, |lb|) \cdot N_L \cdot \sum_{j=0}^{N_S} \delta_j \right) \notag \\
\text{s.t.} \quad 
& \begin{array}[t]{ll}
lb \cdot \alpha_{ij} \leq \xi_{ij} \leq ub \cdot \alpha_{ij}, & \alpha_{ij} \leq \delta_j \\[1.2mm]
\sum_{i,j=1}^{N_t,N_S} \alpha_{ij} \leq K_\alpha, & \sum_{j=1}^{N_S} \delta_j \leq K_\delta \\[1.2mm]
\xi_{ij} \in \mathbb{R}, & \delta_j, \alpha_{ij} \in \{0,1\}
\end{array}
\end{align}

The optimisation includes binary variables $\alpha_{ij}$, whose value constrains the boundaries for the search of $\xi_{ij}$, and their cumulative contribution is selected through the parameter $K_\alpha$. Additionally, the need to activate only a fixed number of columns within the matrix $\boldsymbol{\Xi}$ requires the inclusion of a further auxiliary binary variable $\delta_j$ that constraints the value of $\alpha_{ij}$ to be active only if a given column is selected (i.e., $ \alpha_{ij} \leq \delta_j $).

\subsection{MILP reformulation of the problem}
Including the absolute error rather than the squared error for evaluating the model performance, as the classical SINDy algorithm, allows for the problem reformulation into an MILP configuration. The problem described in \eqref{eq:7} can be converted into a MILP problem to increase its treatability through the inclusion of auxiliary variables $y$ and $z$ representing the absolute value of the errors and the model coefficients, respectively. Therefore, constraints on the values of $y$ and $z$ are applied as follows $y_{ij} \geq (h_{exp} - X_L \cdot \boldsymbol{\Xi})_{ij}$ and $y_{ij} \geq -(h_{exp} - X_L \cdot \boldsymbol{\Xi})_{ij}$, $z_{ij} \geq \xi_{ij}$ and $z_{ij} \geq -\xi_{ij}$. 
Additionally, the presence of the regularisation term on the number of DDM equations (i.e., $\lambda_{1,\xi} \cdot \max(|ub|, |lb|) \cdot N_L \cdot \sum_{j=0}^{N_S} \delta_j$) makes the optimisation function non-smooth and increases the resolution complexity. Therefore, the optimisation is formulated including a further discrete positive optimisation variable $s$ to improve the problem treatability; its value is constrained between $0$ and $N_s$ and satisfies $\sum_{i=0}^{N_S} \delta_i-s \leq 0$. Willis and von Stosch \cite{Willis_etal_2017} employed this formulation for model identification tasks, while other authors successfully used it for process control tasks \citep{DeOliveira_etal_1994, Richards_2015}.
After applying all the abovementioned transformations, the problem is reformulated as a MILP optimisation in the variables $\boldsymbol{\Xi}, \boldsymbol{A}, \boldsymbol{\Delta}, \boldsymbol{Y}, \boldsymbol{Z}, s$ upon the selection of hyperparameters $\lambda_{1,\xi}$, $K_\alpha$ and $K_\delta$, with a system having $N_S$ states and $N_L$ functions contained in the function library $\boldsymbol{X}_L$. Therefore, the problem is reformulated as
\begin{align}
\label{eq:8}
\min_{\boldsymbol{\Xi}, \boldsymbol{A}, \boldsymbol{\Delta}, \boldsymbol{Y}, \boldsymbol{Z}, s} \quad
& \sum_{i,j=0}^{N_t,N_S} y_{ij} 
+ \lambda_{1,\xi} \left( \sum_{i,j=0}^{N_L,N_S} z_{ij} + \max(|ub|, |lb|) \cdot N_L \cdot s \right) \notag \\
\text{s.t.} \quad
& 
\begin{array}[t]{ll}
y_{ij} \geq (h_{exp} - X_L \cdot \boldsymbol{\Xi})_{ij}, & y_{ij} \geq -(h_{exp} - X_L \cdot \boldsymbol{\Xi})_{ij} \\[1.2mm]
z_{ij} \geq \xi_{ij}, & z_{ij} \geq -\xi_{ij} \\[1.2mm]
lb \cdot \alpha_{ij} \leq \xi_{ij} \leq ub \cdot \alpha_{ij}, & \alpha_{ij} \leq \delta_j \\[1.2mm]
\sum_{i=0}^{N_S} \delta_i - s \leq 0, & \sum_{i,j=1}^{N_t,N_S} \alpha_{ij} \leq K_\alpha \\[1.2mm]
\sum_{j=1}^{N_S} \delta_j \leq K_\delta, & y_{ij}, z_{ij} \in \mathbb{R}^+ \\[1.2mm]
\xi_{ij} \in \mathbb{R}, & s \in \{0, \dots, N_s\} \subset \mathbb{Z}_+ \\[1.2mm]
\multicolumn{2}{l}{\delta_j, \alpha_{ij} \in \{0,1\}}
\end{array}
\end{align}
where $\boldsymbol{Y}$ and $\boldsymbol{Z}$ are the matrices containing the slack variables $y$ and $z$.
Often, ODE comes with significantly different value scales for the various system states, the magnitude of which can be further amplified when computing the library function for the generation of the $X_L$ matrix. This can potentially give bias toward the identification and correction of states with a bigger magnitude, as well as bias on the scale of the identified parameter value of the DDM. As a matter of fact, large $\xi$ parameter values could be requested on library functions with small values. However, this is penalised by the L1 regularisation of the formulation in \eqref{eq:7} and \eqref{eq:8}. Therefore, to avoid adding bias to the search and favour equal parameter identification, the values of the $X_L$ columns were normalised based on the maximum of their absolute values
\begin{equation}\label{eq:9}
X_{L,ij}^{sc} = \frac{X_{L,ij}}{\max(|X_{L,j}|)}
\end{equation}
the $X_L^{sc}$ values were used in \eqref{eq:8} in place of $X_L$.

\section{Experimental validation}
\subsection{Investigation cases and methodology description}
\label{sec:investigation_case_description}
The optimisation problem \eqref{eq:8} is applied to three dynamic cases obtained from the literature to verify the capabilities of the SINDybrid methodology, including a Meerwein arylation reaction executed in catalytic batch conditions, a bioethanol production fermentation process using Saccharomyces cerevisiae in continuous conditions and a Lotka-Volterra oscillator.

The \textit{Meerwein arylation reaction} case describes the temporal evolution of 2 reactant concentrations (i.e., C\textsubscript{A} and C\textsubscript{B}), 1 desired product concentration (i.e., C\textsubscript{P}) and 1 side-product concentration (i.e., C\textsubscript{S}), described by a set of 4 ODEs. Besides the 4 concentration states, the system includes two further experimental run conditions, which remained constant during the execution of the experiment, namely the experimental temperature (T) and the concentration of catalyst (C\textsubscript{cat}). The kinetic parameters and the running conditions boundaries for the Meerwein arylation reaction case were obtained from the literature \cite{Shukla_Atapalkar_Kulkarni_2020}. The \textit{continuous fermentation process} describes the temporal evolution of a biomass concentration, a substrate concentration and a product concentration, namely X, S and P. The dynamic was described by 3 ODEs, including biomass growth inhibition caused by product and biomass accumulation. Further, the model considers 2 experimental run conditions to characterise the system dynamics, namely the dilution factor (D) and the inlet concentration of substrate (S\textsubscript{A}). The parameters, the models, and the running condition boundaries were obtained from the literature \cite{Elmer_etal_2013}. The \textit{Lotka-Volterra oscillator} described the temporal evolution of a prey-predator system, with a set of 2 ODEs; all the solutions of this system are periodic with two steady equilibrium points \cite{Wangersky_1978}. Four parameters characterise the system dynamics; the validation cases reported in this work were all set to 1. The symbolic representation of the employed models and further information about the initial conditions are reported in the equations \eqref{eq:s1}, \eqref{eq:s3}, and \eqref{eq:s5} within the SI section of this work.

All validation optimisations were conducted with relaxed constraints on the DDM cardinality and the number of DDMs included (i.e., with \(K_\alpha = \infty\) and \(K_\delta = \infty\)). During the validation procedure, a known deviation was applied to represent the epistemic uncertainty of the system using a specified function for each state within the system. Following this, the optimisation process was executed, and the results were evaluated based on several factors: the number of activated DDMs, their position in the FPM, and the quality of the approximation for the epistemic uncertainty. For detailed information on the deviations applied to the model equations, please refer to the SI section of this work, particularly in equations \eqref{eq:s2}, \eqref{eq:s4}, and \eqref{eq:s6}.
The investigation aims to verify the robustness of the methodology in terms of experimental noise level, the number of experiments employed for the identification, and the number of time samples in each experiment. The nominal validation case employed 6 experiments for the HM identification and training, and 2 for the evaluation of the generalisation capabilities of the model. The training and the testing points are obtained using a random sampling strategy. The initial conditions for the ODE were sampled using Latin Hypercube sampling, which involved both experiments in the training set and testing set. For the cases containing experimental run conditions (i.e., Meerwein arylation reaction and the continuous fermentation process), the points were sampled using a full-factorial approach within the boundaries suggested by the literature, and through random selection of the conditions to match the number of given experiments. Further information about the running condition management for each case is reported in the SI section of the work.
Each nominal experiment was executed with a uniform noise on the system states of 10\%, and the algorithm was executed with a regularisation factor $\lambda_{1,\xi}=3$. Nevertheless, an investigation into the effect of these values was executed and can be found in the SI section of this work, in Figure \ref{fig:4s}, Figure \ref{fig:5s}, and Figure \ref{fig:6s}. The quality of the SINDybrid approach was investigated through 3 main parameters, namely 1) the ability for the algorithm to identify the correct position of the epistemic uncertainty evaluated comparing the actual epistemic uncertainty position and the identified one, 2) the quality of the model to capture the information contained in the epistemic uncertainty evaluated calculating the coefficient of determination ($R^2$) metric on the HM prediction and 3) the accuracy of the overall model in quantitative terms evaluated through the mean absolute error (MAE) of the developed HM.

The effect of the data quality and quantity on the HM identification was executed through investigating: 1) the effect of the noise on the experimental data, 2) the effect of the number of batches included in the training set, and 3) the effect of the time-sampling rate in the experimental data. The first influences the quality of the training data; increased noise compromises the integrity of the data from which the algorithm identifies the DDMs and affects the calculation of the FPM in \eqref{eq:2}, reducing identification capability and overall model accuracy. The second and third factors influence the quantitative information contained in the matrix $h_{exp}$ by reducing the number of rows available for evaluating the obtained DDMs in the identification phase. Specifically, the number of batches impacts the diversity of experiments and conditions observed during optimisation, whereas the time-sampling rate contributes to the sparse representation of information within the system.

\subsection{Results of the investigation}
\subsubsection{Effect of the experimental noise}
\begin{figure}[H]
    \centering
    \includegraphics[width=\textwidth]{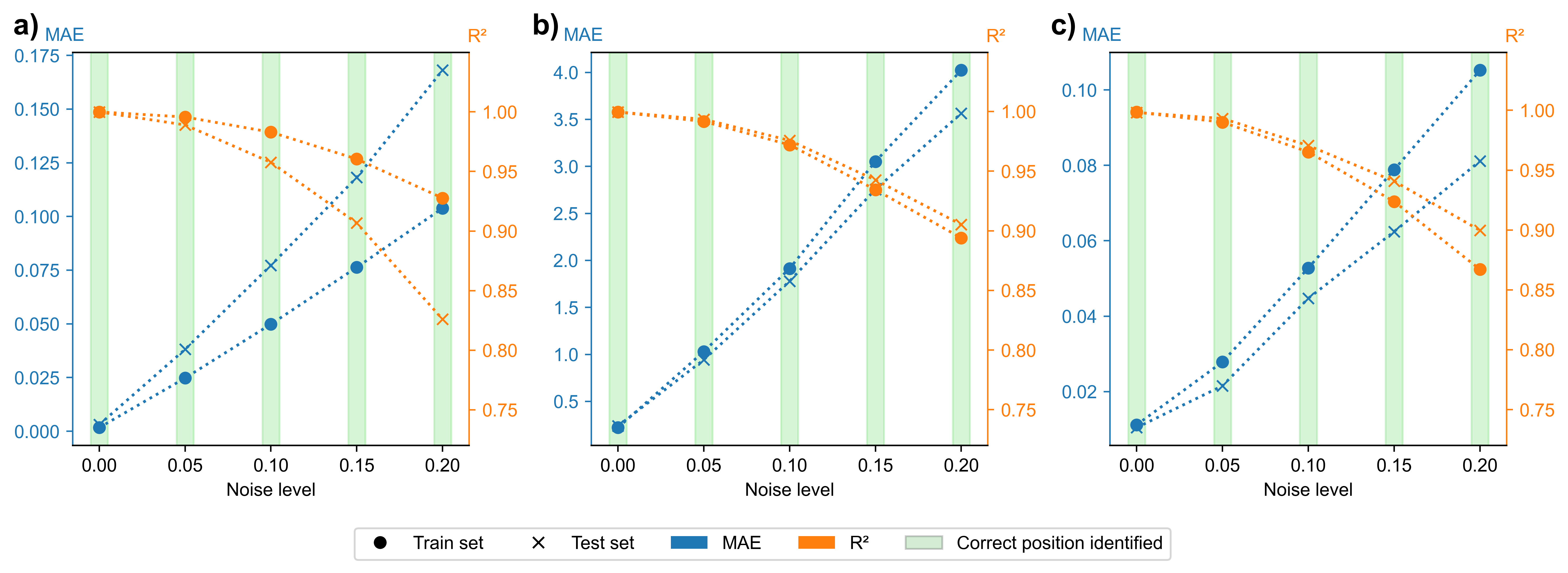}
    \caption{Effect of the experimental noise on the HM identification and training for the a) Meerwein arylation, b) Continuous fermentation and c) Lotka-Volterra oscillator cases. The points are an average of the metrics on all the process variables for all the deviations applied to the model. The vertical bends report the success of the identification; if any of the identification fails on any level, the bend is red. All the identification results were stable both on the train and test sets until an experimental noise of 15\%, with great model accuracy. The accuracy decreases between the training and testing sets, and is more marked in the Meerwein arylation compared to the other cases. In contrast, the metrics in the other two cases remain similar between the training and testing sets.}
    \label{fig:1}
\end{figure}

Figure \ref{fig:1} depicts the effect of experimental noise on HM identification and training across the three validation cases. For each case, perturbations were individually applied to all model equations, and the corresponding SINDybrid identification errors were recorded. In Figure \ref{fig:1}, each panel corresponds to one of the three test cases. Each point represents the error obtained across multiple runs at a given noise level. Each run sets the synthetic epistemic uncertainty location on a different model equation to generate epistemic uncertainty. The green vertical bands indicate that the algorithm successfully identified all applied deviations; identification was considered successful when the algorithm correctly identified the position of all the applied deviations for a given case. The accuracy metrics (i.e., $R^2$ and $MAE$) were computed for all variables across the validation cases at each noise level and averaged. In Figure \ref{fig:1}, one can see how the epistemic uncertainty location identification was successful and robust across the various tested scenarios at all the noise levels. Additional tests conducted with noise levels up to 35\% revealed only one case of failed identification, namely the continuous fermentation case, where performance deteriorated beyond 20\% noise. Detailed representations of these experiments can be found in the SI, Figures \ref{fig:1s} to \ref{fig:3s}. The Meerwein arylation and Lotka–Volterra oscillator hybrid models were successfully identified and trained, with model accuracy decreasing as the data noise increased. In the SINDybrid framework, experimental noise has a dual influence: 1) it distorts the numerical estimation of time derivatives ${d\boldsymbol{x}}/{dt}$, and 2) it affects the evaluation of $\boldsymbol{f}(\boldsymbol{x}_{exp}, \boldsymbol{r}_{exp})$, both of which are used in computing the matrix $\boldsymbol{h}_{exp}$ used to identify the sparse coefficient matrix $\boldsymbol{\Xi}$. The higher sensitivity of the continuous fermentation case to noise is likely attributable to the presence of exponential and rational functions in the mechanistic model (e.g., $e^X$ and $\frac{X}{K + X}$). These nonlinearities amplify noise-induced errors, leading to substantial deviations from the noiseless model prediction. Nonetheless, the methodology exhibited considerable robustness up to 20\% noise, indicating tolerance to a certain level of aleatoric uncertainty. In contrast, the Meerwein arylation and Lotka–Volterra oscillator models include terms involving the product between system parameters and variables (e.g., $k_1C_A$ and $k_1C_A C_B$), which did not affect the success of model identification but did influence model accuracy. This may be due to the absence of complex nonlinearities that significantly amplify noise in the mechanistic formulation.

The values of the metrics in Figure \ref{fig:1} highlight the high accuracy of the model both on the training set and on the testing sets for experimental noise up to 20\%. The $R^2$ values are above 0.8 for all the test cases, demonstrating how the model captures the trend of the epistemic uncertainty from the data while retaining the information accuracy and transferring it to the testing set. The plots in Figure \ref{fig:1} show an aggregated trend for all the deviations included in the models. However, the trends are the same for all the variables (c.f., Figure \ref{fig:1s} to \ref{fig:3s} in the SI section). The MAE differ quite significantly between the continuous fermentation (i.e., ranging between 0.25 and 4.0) and the other validation cases(i.e., ranging between 0.01 and 0.2), and this is related to the magnitude of the system states; the continuous fermentation has state values in the order of magnitude of 100, while the Meerwein arylation and the Lotka-Volterra oscillator report states with values around 2, with an MAE being in the range 5-10\% of the maximum range reached by the variables.
In Figure \ref{fig:1}, all the curves have comparable metrics on training and testing sets at experimental noise equal to 0, highlighting the capabilities of the algorithm to identify the correct DDM position and generalise to new experiments in these conditions. Increasing the experimental noise caused the model accuracy to decline on both the training and test sets. This outcome is expected, as higher noise levels reduce the effective information content of the data, making it more challenging to extract meaningful patterns. The metrics on the training and testing sets have similar trends for both continuous fermentation and Loka-Volterra oscillator cases (Figure \ref{fig:1}b and Figure \ref{fig:1}c). In contrast, in the Meerwein arylation (Figure \ref{fig:1}a), the error increase is larger for the testing set than the training one, with a significant leap between the two lines starting from an experimental noise of 15\%, showing poor generalisation capabilities of the obtained model. This behaviour could be related to the fact that one of the experiments contained in the testing set of the Meerwein arylation case is outside the model training boundaries with the DDM forced to work in extrapolation from the training set; with small experimental noise within the training set, the algorithm is able to identify the correct deviation equation and makes it flexible on all the generalisation cases. When the experimental error on the training set increases, the solver experiences increased difficulty in identifying the model describing the epistemic uncertainty, generating an approximation function that is less accurate as the noise intensity increases. Such an approximation decreases its quality the higher the experimental error gets, decreasing the extrapolation capabilities of the obtained model and returning the increased deviation in the metrics between the training and testing sets.

\subsubsection{Effect of the experimental run amount}
\begin{figure}[H]
    \centering
    \includegraphics[width=\textwidth]{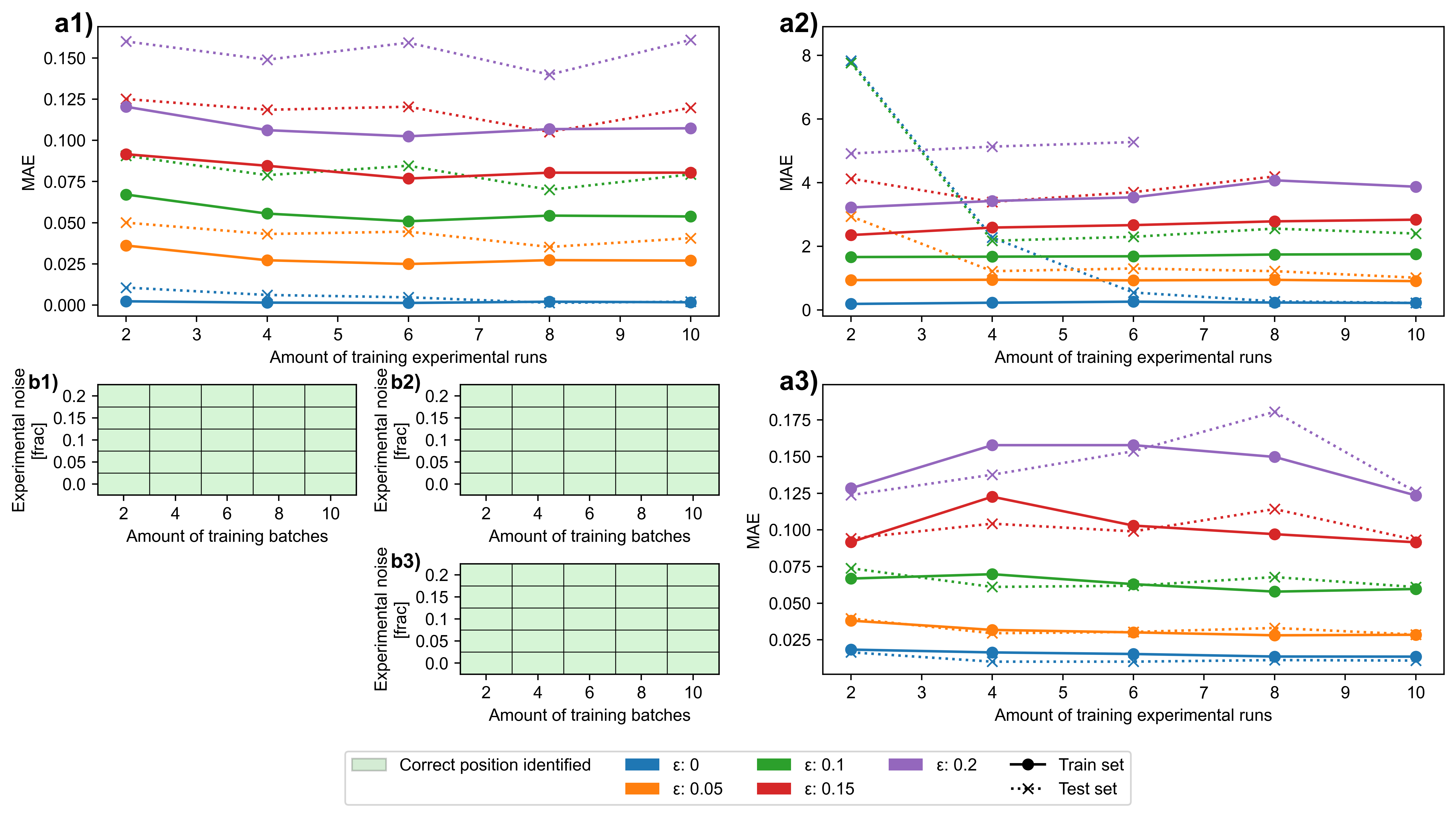}
    \caption{Effect of the number of experimental runs included in the training set on the identification results, for the 1) Meerwein arylation, 2) Continuous fermentation and 3) Lotka-Volterra oscillator cases. Plots a), b) and c) report the trend of the MAE metric with the number of points in the training set, parametric with the experimental error. The points are an average of the metrics on all the process variables for all the deviations applied to the model. Plots a1), b1), and c1) report the success of the position identification. The cells report the success of the identification, taking into account the worst-case scenario for any of the deviations.}
    \label{fig:2}
\end{figure}

Figure \ref{fig:2} depicts the effect of the training batch amount on the epistemic uncertainty identification and DDM training. Each point in Figure \ref{fig:2}a1, \ref{fig:2}a2, and \ref{fig:2}a3 represents the error obtained across multiple runs at a given noise level. Each run sets the synthetic epistemic uncertainty location on a different model equation to generate epistemic uncertainty. The metrics were computed for all variables across the validation cases at each noise level and averaged. The green squares in the grid of Figure \ref{fig:2}a1, \ref{fig:2}a2, and \ref{fig:2}a3 depict the successful location identification at various experimental noise and amounts of training batches. For this analysis, 12 batches were simulated, with the training batches randomly selected among these; the remaining experiments constituted the test set and were used to assess the predictive capability of the model. The figure reports the MAE values following the ODEs integration varying the number of batch runs used for the identification (Figures \ref{fig:2}a1, \ref{fig:2}a2 and \ref{fig:2}a3), and the success of the identification for the epistemic uncertainty (Figures \ref{fig:2}b1, \ref{fig:2}b2 and \ref{fig:2}b3) for all the validation cases employed in this study. In Figures \ref{fig:2}a1, \ref{fig:2}a2 and \ref{fig:2}a3, increasing the number of experiments in the training set has little effect on increasing the model accuracy, meaning that the information retrieved from the system is invariant to the number of experiments employed for the identification. In fact, with the increase in the number of experiments, the metrics do not show any significant improvement in the training set and the testing set. The invariance of the obtained information with the number of experimental batches does not change when the experimental noise is observed for all the cases on the training and testing sets. 

The exception to this behaviour is the continuous fermentation process at low experimental noise (i.e., 0\% and 5\%), as shown in Figure \ref{fig:2}a2. For these cases, the MAE on the training set does not change when varying the amount of employed experimental data; on the other hand, the MAE on the testing set significantly decreases. Such behaviour highlights the model overfitting on the training set for a low number of experiments. It could be related to the rational functions included in the continuous fermentation process that increase the sensitivity of the model state variation and, therefore, increase the possibility of model overfitting. For this reason, if strong non-linearities are present in the FPM, the algorithm requires more experiments to improve the generalisation capabilities of the obtained model. 
Additionally, in the continuous fermentation case, at high experimental noise level (i.e. 15\% and 20\%) and high amount of experiments in the training set (i.e., 8 and 10 experiments in the training set) the identification of the epistemic uncertainty location succeed (Figure \ref{fig:2}b2) but the ODE integration on the testing set fails because of numerical integration problems. Also, for this behaviour, we hypothesise that it is an effect of rational functions in the FPM and the large number of experiments used in the training set. In fact, by increasing the experimental noise, the error when computing the FPM values increases compared to the other cases because of some variables located in the denominator. Increasing the FPM computation error affects the calculation of the $\boldsymbol{h}_{exp}$ matrix that is employed in the optimisation function. Additionally, the higher number of experimental points increases the optimisation complexity because of the high number of cases to minimise. Such a dual effect results in the model overfitting on the training set, which results in an accurate resolution on such a set, but diminishes the generalisation capabilities of the testing set, resulting in unstable ODE integrations. As a matter of fact, such behaviour does not show when the noise level is lower (i.e., lower error propagation in the FPM calculation) and when the amount of data points is decreased (i.e., the number of cases considered in the optimisation is lower). Moreover, Figure \ref{fig:2}a2 depicts how the model fails the ODE resolution on the test set at 8 training experiments when the experimental noise is 15\% while this value decreases to 6 when the experimental noise is at 20\%, proving the presence of the combination of the two effects. Therefore, de-noising the profile when working with FPM with rational functions is strongly advised to increase the generalisation capabilities and the ODE integration success of the obtained HM. Despite the abovementioned effects, the algorithm was able to identify the location of the epistemic uncertainty within the model for all the cases employed in this work, also for the ones in the continuous fermentation case where the ODE integration on the testing set was failing, as reported in Figures \ref{fig:2}b1, \ref{fig:2}b2 and \ref{fig:2}b3. Therefore, in the best-case scenario, SINDybrid returns a model able to generalise on novel cases, while in the worst-case scenario, it provides the location of the epistemic uncertainty, being a valuable tool for further investigating the system and using alternative techniques for the HM training.

\subsubsection{Effect of time-samples amount}
\begin{figure}[H]
    \centering
    \includegraphics[width=\textwidth]{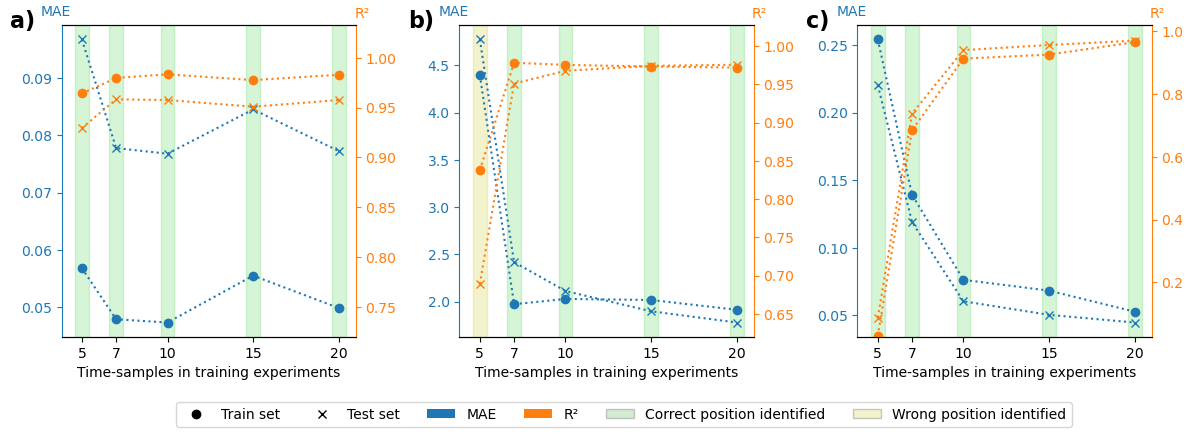}
    \caption{Effect of the number of time samples per experiment for the runs included in the training set for the a) Meerwein arylation, b) Continuous fermentation and c) Lotka-Volterra oscillator cases. The points are an average of the metrics on all the process variables for all the deviations applied to the model. The vertical bends report the success of the identification, if any of the identification failed on any of the levels, the bend is red. For the Meerwein arylation reaction and the continuous fermentation, the model accuracy reaches a plateau after 7 time points, while this value increases to 10 for the Lotka-Volterra oscillator. Additionally, all the positions of the epistemic uncertainty are correct in all the cases for all the time samples except for the continuous fermentation case, which fails to identify with 5 data points.}
    \label{fig:3}
\end{figure}

Figure \ref{fig:3} reports the identification and training performance of SINDybrid varying the time samples within the experimental set. The metrics were calculated for all the system states in a validation case (i.e., set of deviations) and averaged. The identification was considered successful in case the location of the epistemic uncertainty model was correctly identified by the algorithm for all the applied deviations. The individual identification performance for each of the tested epistemic uncertainty positions and validation cases is reported in the SI of this work in Figure \ref{fig:7s} to Figure \ref{fig:9s}. In the Meerwein arylation and Lotka-Volterra oscillator cases, the identification of the epistemic uncertainty position was successful for experiments with only 5 time samples within the training set. On the other hand, in the continuous fermentation, with 5 time samples, a wrong position was identified for one of the cases. More specifically, this happened when the deviation mimicking the epistemic uncertainty source was placed on the biomass change (i.e., ${dX}/{dt}$). Here, the algorithm identified that the epistemic uncertainty was located in the equation describing the change of the product (i.e., ${dP}/{dt}$). Such behaviour is related to the low number of data points employed to compute the identification matrix $\boldsymbol{h}_{exp}$ and the structure of the FPM describing the continuous fermentation case. The dynamics of $X$ and $P$ are highly correlated since the biomass concentration is included in the production of $P$ through $r_P$, and the product concentration is included in the $r_x$ equation as the poisoning effect. With this high correlation and with a small number of data points to differentiate between the two effects, SINDybrid assigns the deviation to the condition that corrects the most the deviation distribution on the data (c.f., Section \ref{sec:why_does_the_algorithm_work}). This behaviour was also observed for systems without experimental noise, which eliminates the hypothesis of an effect of the noise. Increasing the number of time sample points to 7 in the continuous fermentation case, the SINDybrid approach correctly identifies the epistemic uncertainty position.
In Figure \ref{fig:3}, it can be observed that the accuracy of the model obtained with the SINDybrid approach increases as the number of data points increases. This is expected since adding new data points brings more examples in the matrix $\boldsymbol{h}_{exp}$ that are used by the optimisation to fine-tune the model. Yet, the increase in model accuracy is a different order of magnitude for the investigated cases. On the Meerwein arylation and the continuous fermentation cases, the improvement is quite marginal, ranging from a 10-20\% decrease of the MAE and with substantially unchanged $R^2$ values. On the other hand, in the continuous fermentation case, one can observe a marginal improvement of the metrics from experiments with more than 5 data points. At 5 time-samples, SINDyrid fails the epistemic position identification, which lets the metrics significantly worsen both on the training set and on the testing set compared to the cases with a higher number of time-points. In contrast to the Meerwein arylation and continuous fermentation cases, the Lotka–Volterra oscillator shows significant improvement when moving from 5 datapoints to 10 datapoints with an MAE decrease of around 70\%, being the case that benefits the most, increasing the time samples in the experiments. This is related to the periodic nature of the system; it shows a lower co-domain variance in the output space, allowing a higher amount of possible solutions fitting the available points compared to non-periodic cases (i.e., Meerwein arylation and continuous fermentation). Therefore, adding an additional time-point to the system decreases the number of possible solutions until the algorithm is able to reach the optimal solution from experiments with more than 10 time samples. This phenomenon is not observed in Meerwein arylation and continuous fermentation because there is a higher co-domain variance of the data points and a lower number of possible solutions. There, adding new data points favours the fine-tuning of the parameters around the previously obtained value rather than the identification of a new function. Additionally, the periodic nature of the system restricts the function co-domain variance, decreasing the quality of the available information; therefore, adding a novel time point brings more information in a periodic case than a non-periodic one. Therefore, when using the SINDybrid approach, we suggest increasing the number of time samples in the case of periodic FPM. Nevertheless, the model accuracy reaches a plateau for all the validation cases for the number of time samples equal to 10. Moreover, from Figure \ref{fig:3}, the model accuracy on the training and the testing sets is comparable, with error metrics very close to each other within the entire time-samples range, highlighting the great generalisation capabilities of the obtained model for all the time-samples investigated in this study.

\subsubsection{Quantitative experiments}
\begin{figure}[H]
    \centering
    \includegraphics[width=\textwidth]{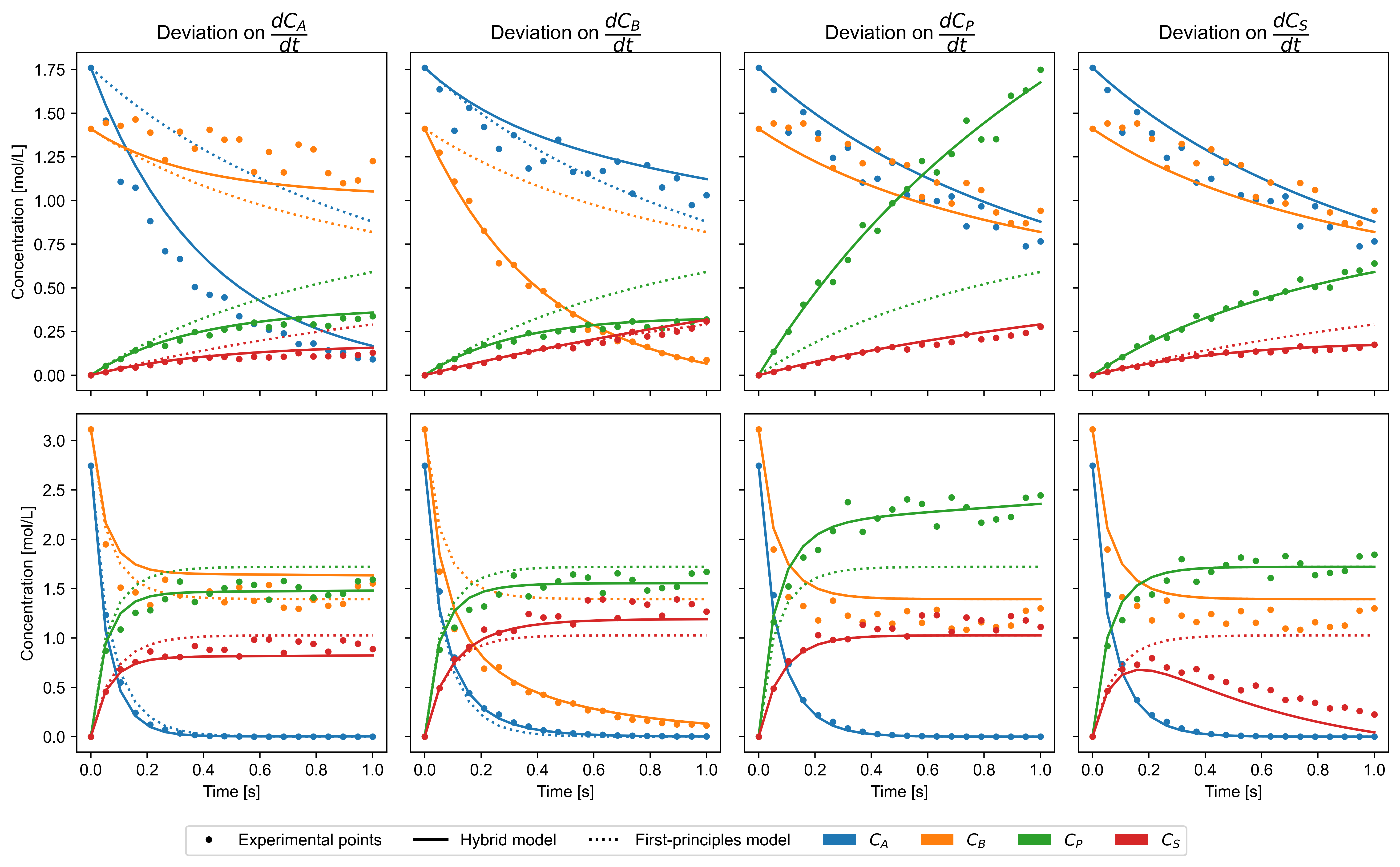}
    \caption{Quantitative model prediction on testing points for the Meerwein arylation experiments for the various deviations happening in the system. Dots are experimental data, the dotted lines are the FPM predictions, and the continuous lines are the predictions of the obtained hybrid model.}
    \label{fig:4}
\end{figure}

\begin{figure}[H]
    \centering
    \includegraphics[width=0.8\textwidth]{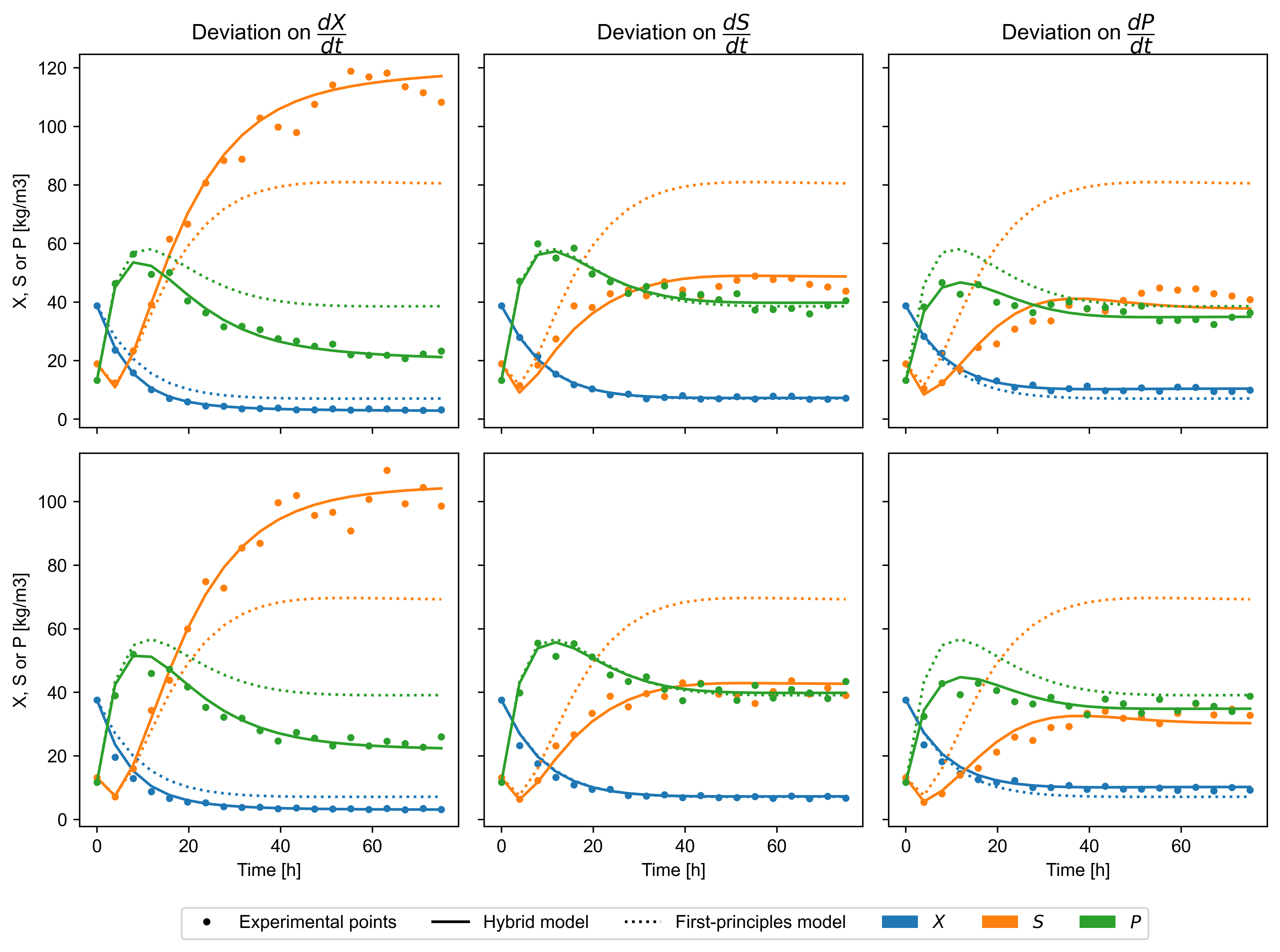}
    \caption{Quantitative model prediction on testing points for the continuous fermentation experiments for the various deviations happening in the system.  Dots are experimental data, the dotted lines are the FPM predictions, and the continuous lines are the predictions of the obtained hybrid model.}
    \label{fig:5}
\end{figure}

\begin{figure}[H]
    \centering
    \includegraphics[width=0.7\textwidth]{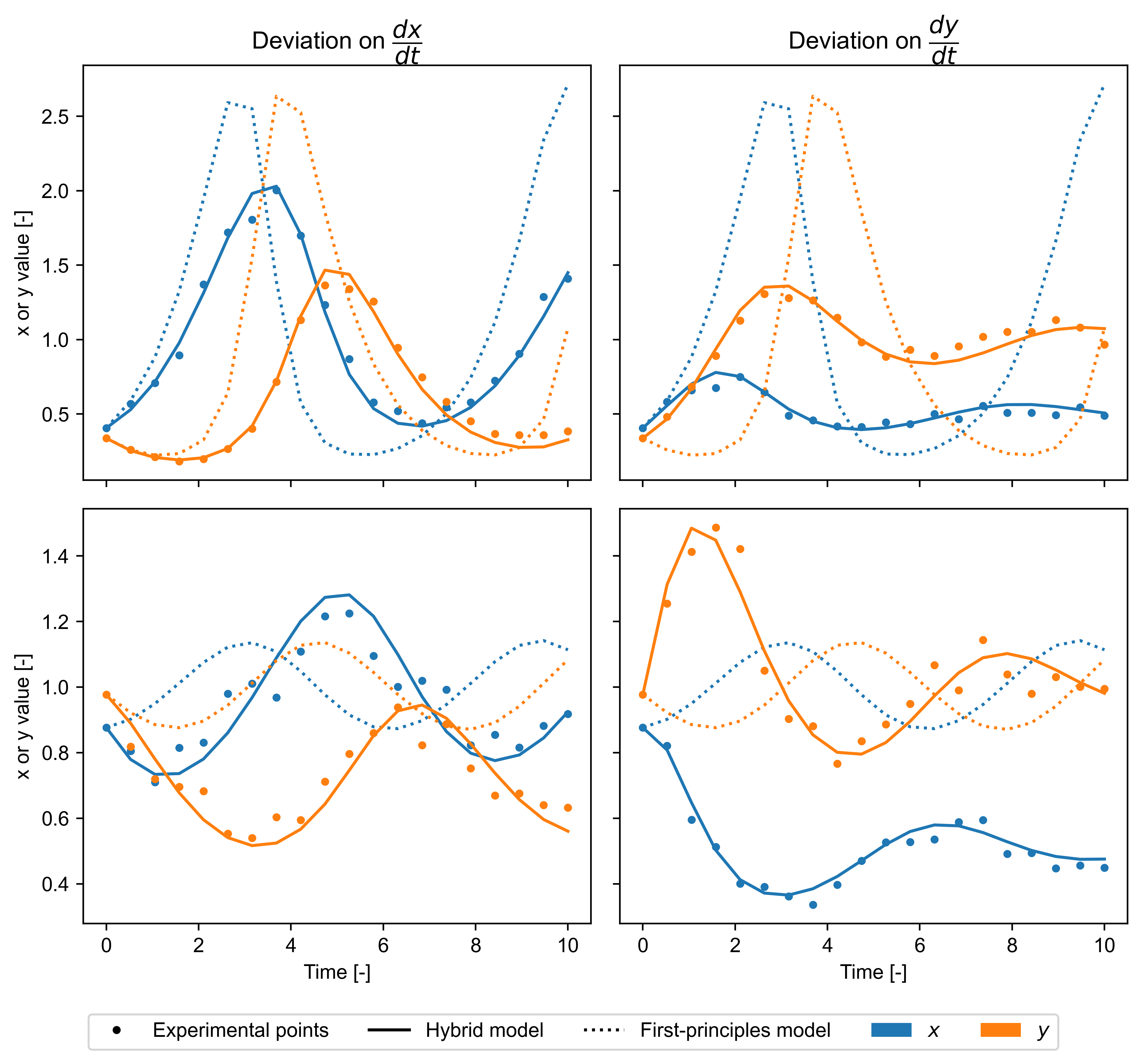}
    \caption{Quantitative model prediction on testing points for the Lotka–Volterra oscillator experiments for the various deviations happening on the system.  Dots are experimental data, the dotted lines are the FPM predictions, and the continuous lines are the predictions of the obtained hybrid model.}
    \label{fig:6}
\end{figure}

Figure \ref{fig:4}, Figure \ref{fig:5}, and Figure \ref{fig:6} depict the quantitative performance of HM obtained through the SINDybrid approach for the various validation cases employed in this study. The figures depict the results for the Meerwein arylation case, the continuous fermentation case, and the Lotka–Volterra oscillator, respectively. The experiments reported in these figures come from the testing set at the nominal conditions (c.f., Section \ref{sec:investigation_case_description}), reporting the model prediction capabilities. For all cases and applied deviations, predictions of the HM overlap experimental points, regardless of the epistemic uncertainty position. All generated HM display significant deviations from the FPMs results while demonstrating high accuracy in predicting experimental data across all positions. However, not all the variables predicted by the HM deviate with the same magnitude from the FPM; this mainly depends on how the investigated state and the one that is deviated by the DDM are correlated. More specifically, if the deviation is applied to a state with higher interactions with the other system states, the majority of the states in the system will deviate from the FPM dynamics. For instance, this is the case of $C_A$ in the Meerwein arylation or $X$ in the continuous fermentation case. Therefore, the experimental data and the FPM predictions will differ on all the affected variables. On the contrary, if the epistemic uncertainty is applied to a state that is weakly linked to other system states, only part of the variables deviate from the FPM predictions. In this study, this phenomenon was only observed for the variables $C_P$ and $C_S$ in the Meerwein arylation (Figure \ref{fig:4}). This non-homogeneous behaviour within the system gives a further layer of correlation between the variables, increasing the complexity of the optimisation. However, SINDybrid overcomes the issue related to correlated dynamics since it works in the gradient domain. Such a solution translates to highly accurate models independently of the number of correlation variables.

In this regard, one can evaluate this effect in Figure \ref{fig:4} comparing the identification performed on the cases where $C_A$ or $C_B$ was deviating to the ones where $C_P$ or $C_S$ was deviating. In the first case, the deviated states were highly connected to the others since they describe the dynamics of the reactants, whereas in the second case, there was no connection between the variables because of their role as products of the reaction system. Despite this fundamental difference, the methodology was able to identify the right position of the epistemic uncertainty and compensate for it, achieving an average $R^2$ of 0.928 for the first two cases and 0.894 for the last two. The same applies to the continuous fermentation process prediction reported in Figure \ref{fig:5}. In this case, the dynamics of $X$ and $P$ are highly interacting with all the variables in the system because of their effect on the biomass growth rate, while the $S$ value is less interacting with the others. However, also in this case, the identification showed outstanding capabilities with an average $R^2$ of 0.948 for the cases deviating on $S$ and $P$ dynamics and an $R^2$ of 0.950 for the $S$ dynamic. For the Lotka–Volterra oscillator, the system states typically interact to produce stable orbits with oscillatory behaviour in the absence of epistemic uncertainty. However, the introduction of epistemic uncertainty disrupts this equilibrium, causing the model to converge toward a stable point. The HM generated using the proposed algorithm accurately captures and compensates for the deviation, significantly adjusting the model dynamics to match the experimental data.

Overall, the abovementioned validation cases highlight the robust generalisation capabilities of the HMs derived using SINDybrid techniques. They not only predict system behaviour accurately but also correctly identify the position of epistemic uncertainty, thereby supporting further model investigation and improved experimental design.

\subsection{Why does the algorithm work?}
\label{sec:why_does_the_algorithm_work}
\begin{figure}[H]
    \centering
    \includegraphics[width=0.7\textwidth]{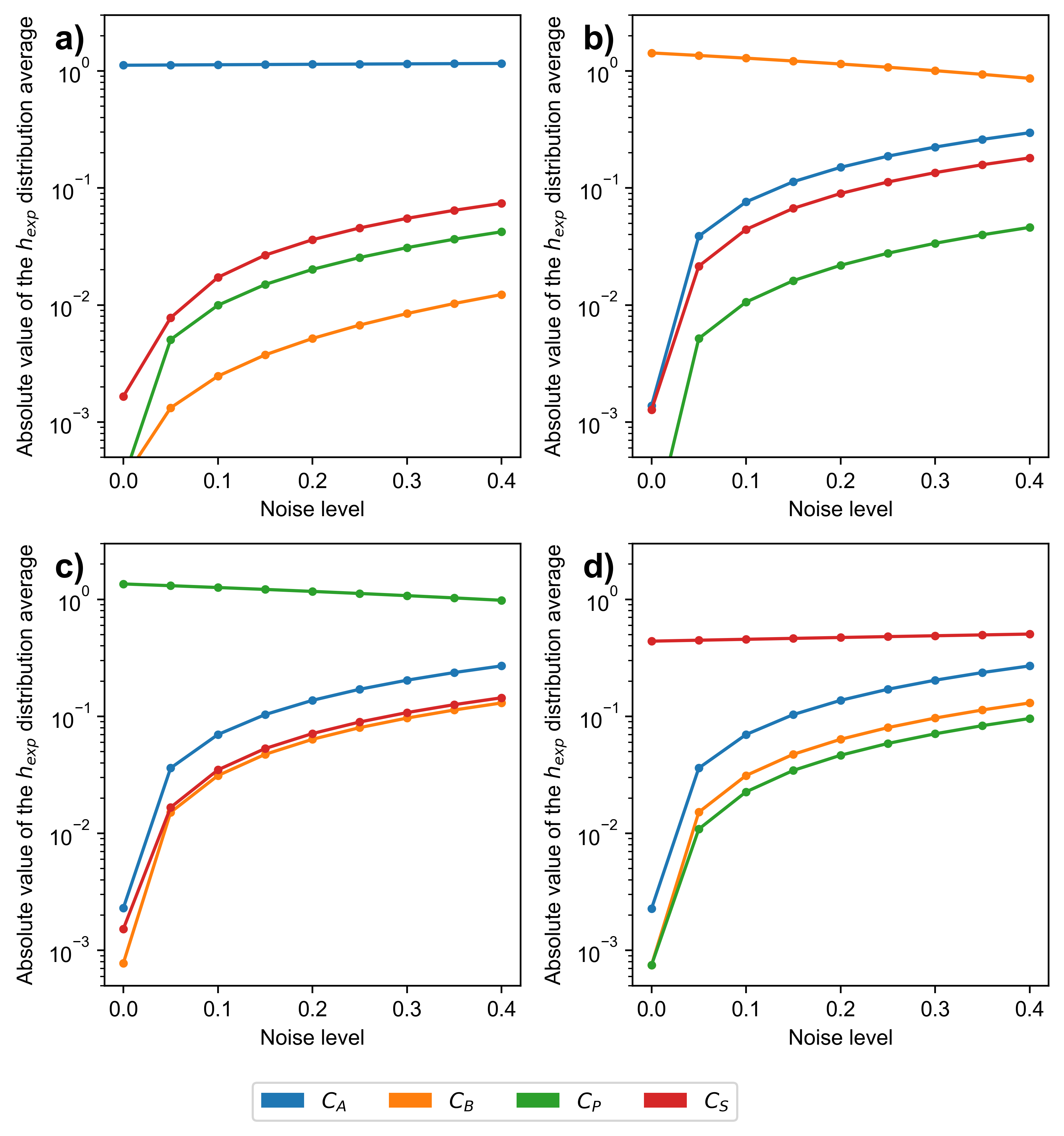}
    \caption{Values of the distribution average for the system deviation in the $h$ matrix at various noise levels for the Meerwein arylation case. a) Epistemic uncertainty on $C_A$ dynamic, b) Epistemic uncertainty on $C_B$ dynamic, c) Epistemic uncertainty on $C_P$ dynamic, d) Epistemic uncertainty on $C_S$ dynamic.}
    \label{fig:7}
\end{figure}

To understand how the algorithm performs the identification, one should observe $\boldsymbol{h}_{exp}$ matrix described in \eqref{eq:2}. In fact, the algorithm applies a SINDy approach, selecting the variable that shows the major deviation from the Gaussian distribution. Looking at \eqref{eq:2}, a dynamics not affected by any epistemic uncertainty should be, in principle, represented by zero columns in the matrix $\boldsymbol{h}_{exp}$. However, in practice, this value cannot be identically zero because of the numerical error derived from the numerical derivatives calculation using the experimental data and the presence of the noise within the measurements; both these errors can be assumed to be centred on zero because of their random nature. Therefore, in reality, the variables whose dynamic is not deviated by any epistemic uncertainty are distributed around zero with a Gaussian shape. On the other hand, in case a variable dynamic is running under epistemic uncertainty, its distribution is not centred around zero because of the presence of the epistemic uncertainty not taken into account in the FPM, which can be represented as a function of the system states, nulling the hypothesis of aleatoric uncertainty. In this scenario, increasing the time-sampling points in the training set is beneficial both to decrease the error when calculating the numerical derivatives and to have more samples to evaluate the distribution of the $\boldsymbol{h}_{exp}$ columns. The behaviour differentiating a system state running with epistemic uncertainty from one running under aleatoric uncertainty is reported in Figure \ref{fig:7}. Here, the distributions of the deviation dynamics of the system states (i.e., $\boldsymbol{h}_{exp}$) are reported against the experimental noise at every location of the epistemic uncertainty. For brevity and as an example note, Figure \ref{fig:7} depicts only the Meerwein arylation case, but what will be said in this section is valid for all the cases investigated in this work (Figure \ref{fig:10s} and Figure \ref{fig:11s}). Despite the presence of the error related to the numerical derivation, from this figure at a noise level of 0, it is possible to clearly differentiate between variables running under epistemic uncertainty and those that are not. As a matter of fact, for all the applied deviations, the difference in the averages is one order of magnitude larger, making this differentiation very clear. The MILP algorithm benefits from this scenario since this implies that the application of a function on the system state running under epistemic uncertainty would return a major change in the output distribution and, therefore, a major effect on decreasing the optimisation function. Additionally, increasing the experimental noise, the shift between the state with epistemic uncertainty and those that are not is still relevant. This is dictated by the order of magnitude of the experimental noise and the epistemic uncertainty deviation; the epistemic uncertainty deviation is proportional to the value of the system states and has an order of magnitude that is comparable to them, while the noise is typically a fraction of their value (i.e., 0-40\% in the examples reported in Figure 7). For this reason, in Figure \ref{fig:7}, increasing the experimental noise has a negligible effect on the variables running under epistemic uncertainty, while the noise returns a linear effect on the variables that are not affected by the epistemic uncertainty. This behaviour could explain why the algorithm is highly stable with the experimental noise added to the data.

\subsection{Limitation of the proposed approach}
Despite its robustness and ability to produce accurate models, the SINDybrid methodology has several limitations that may restrict its broader application in engineering contexts. As shown in this work, SINDybrid relies on the classical SINDy framework to approximate the deviations contained in the $\boldsymbol{h}_{\text{exp}}$ matrix. Consequently, it inherits the limitations of this approach, most notably, the restriction to function libraries with fixed structures and linear-in-the-parameters formulations.
Moreover, the current implementation addresses only epistemic uncertainty at the equation level without accounting for uncertainties embedded in model parameters. While, in principle, parameter-level deviations could also be addressed using SINDy-based techniques, this aspect was outside the scope of the present study. The algorithm, therefore, operates under the assumption that the FPM parameters have been correctly identified. This limitation can be mitigated by performing a dedicated parameter fitting stage prior to applying SINDybrid, thereby ensuring the best use of the mechanistic information before uncertainty identification begins.
Another constraint of the methodology is its reliance on full-state measurements. Since the algorithm infers deviations based on the observed system dynamics, the absence of measurements for one or more states may affect its ability to localise and compensate for epistemic uncertainty correctly. Nonetheless, if the deviating state is still being measured, the algorithm should, in principle, retain the ability to detect and correct the deviation. This scenario, however, was not explored in this work, which focused exclusively on systems for which all state variables were observed.

\section{Conclusion}
This work presents \textit{SINDybrid}, a novel algorithm for the automatic generation of hybrid models (HMs) in dynamic systems governed by ordinary differential equations (ODEs). The method identifies model equations affected by epistemic uncertainty and compensates for them using a sparse regression framework based on the Sparse Identification of Nonlinear Dynamical Systems (SINDy). By formulating the identification task as a Mixed-Integer Linear Programming (MILP) problem, SINDybrid leverages a constrained optimisation process that incorporates a predefined function library, L1 regularisation, and sparsity constraints to ensure model interpretability and robustness. The algorithm was validated across three representative dynamic systems relevant to chemical engineering: a catalytic batch reaction, a continuous bioethanol fermentation process, and a Lotka–Volterra oscillator. These case studies simulate conditions consistent with real-world experimental constraints in the chemical industry, including limited batch numbers, noisy measurements, and sparse temporal sampling. SINDybrid consistently demonstrated the ability to identify the location of epistemic uncertainty within the model equations, accurately reconstruct the associated deviations, and achieve generalisable models with high predictive performance (test $R^2 > 0.89$ in all cases). Performance analysis revealed that the algorithm remains robust up to 20\% measurement noise, although performance degrades slightly for systems involving exponential or rational dynamics (as in continuous fermentation), where error propagation is more severe. Increasing the number of training batches had a limited impact on accuracy, except in cases with low noise, where overfitting may occur. The identification performance improves with the number of time samples, stabilising at around 10 per experiment. Periodic systems particularly benefit from denser sampling due to a reduced solution space for fitting. Further analysis revealed that equations under epistemic uncertainty exhibit statistically significant deviations from the predictions of the first-principles model. These deviations showed a marked departure from normality, facilitating their detection. Notably, this deviation remained distinguishable even under noise levels up to 40\%, highlighting the robustness of the approach. In summary, SINDybrid provides a general, automated framework for HM construction in dynamic systems, effectively bridging first-principles knowledge with data-driven learning. Identifying and correcting epistemic uncertainty at the equation level reduces modelling complexity while enhancing model fidelity and generalisation. This represents a substantial advancement in the development of automated hybrid models, with direct applicability to chemical engineering systems.

Future work will extend SINDybrid to capture parameter-level uncertainty and integrate recent advances in the SINDy framework, broadening its applicability across a wider range of dynamic modelling scenarios.

\section*{Acknowledgments}
The authors acknowledge funding from the KU Leuven project ``HyPro - Automatic hybrid digital twins for process modelling'' (C3/23/007). The authors declare no conflict of interest.
\section*{Data and code availability}
The data and the code used in this work and supporting its findings are openly available in the KU Leuven GitLab repository \url{https://gitlab.kuleuven.be/cipt/sindybrid}.
\section*{Declaration of generative AI and AI-assisted technologies in the writing process}
During the preparation of this work, the authors used Grammarly to check the text for grammar and ChatGPT to increase the fluency of the text. After using this tool/service, the authors reviewed and edited the content as needed and take full responsibility for the content of the published article.

\bibliographystyle{unsrt}
\bibliography{references}


\clearpage
\resetlinenumber 
\setcounter{page}{1}
\setcounter{section}{0}
\setcounter{figure}{0}
\setcounter{equation}{0}
\renewcommand{\thesection}{s\arabic{section}}
\renewcommand{\thesubsection}{s\arabic{section}.\arabic{subsection}}
\renewcommand{\thefigure}{s\arabic{figure}}
\renewcommand{\theequation}{s\arabic{equation}}
\renewcommand{\thepage}{s\arabic{page}}
\renewcommand{\thefootnote}{\fnsymbol{footnote}}
\begin{center}

  {\LARGE \textbf{SINDybrid: automatic generation of hybrid models for dynamic systems}}\\[1.2em]
  
  {\normalsize
  Ulderico Di Caprio\textsuperscript{1}\footnote{Corresponding author: \texttt{ulderico.dicaprio@kuleuven.be}}, 
  M. Enis Leblebici\textsuperscript{1}\footnote{Corresponding author: \texttt{muminenis.leblebici@kuleuven.be}} \\[0.5em]
  \textsuperscript{1} Centre for Industrial Process Technology, Department of Chemical Engineering, KU Leuven,\\
  Agoralaan Building B, 3590 Diepenbeek, Belgium}
\end{center}

\begin{center}
{\LARGE \textit{Supplementary material}}
\end{center}

\section{Validation cases and deviations}
\subsection{Meerwein arylation}
Equation set  \eqref{eq:s1} reports the first-principle model employed for the Meerwein arylation, while Equation set \eqref{eq:s2} reports the applied deviations to mimic the epistemic uncertainty for the various equations in the first-principle model set.

\begin{equation}
\begin{cases}
\frac{dC_A}{dt} = -k_M \cdot C_A C_B - k_S \cdot C_A \\
\frac{dC_B}{dt} = -k_M \cdot C_A C_B \\
\frac{dC_P}{dt} = k_M \cdot C_A C_B \\
\frac{dC_S}{dt} = k_S \cdot C_A \\
k_M = k_{M,0} \cdot \exp(\frac{198.84}{RT}) \\
k_S = k_{S,0} \cdot \exp(\frac{58.81}{RT}) \\
k_{M,0} = 3.71 \cdot 10^{31} \cdot \exp(0.5498 \cdot C_{cat}) \\
k_{S,0} = 1.06 \cdot 10^8 \cdot \exp(0.7669 \cdot C_{cat})
\end{cases}
\label{eq:s1}
\end{equation}

\begin{equation}
\begin{cases}
\frac{dC_A}{dt}|_{DEV} = -C_A C_B - {C_A}/{C_B} \\
\frac{dC_B}{dt}|_{DEV} = -2 C_B  T^* \\
\frac{dC_P}{dt}|_{DEV} = C_A \cdot(-0.2 C_{cat} - 0.2 T^* + 0.3 C_{cat} T^*) \\
\frac{dC_S}{dt}|_{DEV} = -0.1 C_S \cdot (1 + C_S + C_{cat} + T^* + C_{cat} T^*) \\
T^* = {T}/{273.15}
\end{cases}
\label{eq:s2}
\end{equation}

The initial conditions were sampled using a Latin-hypercube strategy, generating initial conditions for the two reactants, with concentration bounds set between [1, 3] for both. The initial concentrations of the reaction products are fixed at zero. Additionally, the function randomly samples the operating conditions of the system, specifically the catalyst concentration ($C_{cat}$) and temperature ($T$), from the discrete sets \{3, 5, 7\} mol\% and \{313.15, 318.15\} K, respectively.

\subsection{Continous fermentation}
Equation set \eqref{eq:s3} reports the first-principle model employed for the continuous fermentation, while equation set \eqref{eq:s4} reports the applied deviations to mimic the epistemic uncertainty for the various equations in the first-principle model set.

\begin{equation}
\begin{cases}
\frac{dX}{dt} = r_X - D \cdot X \\
\frac{dS}{dt} = D \cdot (S_A - S) - r_S \\
\frac{dP}{dt} = r_P - D \cdot P \\
r_X = \mu_{\text{max}} \cdot \frac{S}{K_S + S} \cdot \exp(K_i \cdot S) \cdot (1 - \frac{X}{X_{\text{MAX}}} )^m \cdot (1 - \frac{P}{P_{\text{MAX}}} )^n \cdot X \\
r_S = \frac{r_X}{Y_X} + \frac{\beta_{ms} \cdot S}{K_{\beta_{s2}} + S} \cdot X \\
r_P = r_X \cdot Y_{P/X} + \frac{\beta_{mp} \cdot S}{K_{\beta_{s1}} + S} \cdot X
\end{cases}
\label{eq:s3}
\end{equation}

\begin{equation}
\begin{cases}
\frac{dX}{dt}|_{DEV} = -0.05 \cdot \frac{X \cdot S}{0.01 + S} \\
\frac{dS}{dt}|_{DEV} = -0.01 \cdot \frac{S}{D} \\
\frac{dP}{dt}|_{DEV} = -0.05 \cdot (X + P) - 5 \cdot 10^{-8} \cdot X \cdot P \cdot S_A
\end{cases}
\label{eq:s4}
\end{equation}

The initial conditions were sampled using a Latin-hypercube strategy with bounds of [30, 50] for the initial concentration of X, [5, 20] for the initial concentration of S, and [9, 15] for the initial concentration of P. Additionally, it randomly samples an equal number of experiments as the Latin-hypercube for the dilution coefficient ($D$) and the inlet concentration of S ($S_A$) with values \{0.08, 0.1, 0.12\} and \{160, 170, 180\} respectively.

\subsection{Lotka-Volterra oscillator}
Equation set \eqref{eq:s5} reports the first-principle model employed for the Lotka-Volterra oscillator, while equation set \eqref{eq:s6} reports the applied deviations to mimic the epistemic uncertainty for the various equations in the first-principle model set.
\begin{equation}
\begin{cases}
\frac{dx}{dt} = (1 - y) \cdot x \\
\frac{dy}{dt} = (x - 1) \cdot y
\end{cases}
\label{eq:s5}
\end{equation}

\begin{equation}
\begin{cases}
\frac{dx}{dt}|_{DEV} = -0.2 \cdot x^2 - 0.1 \cdot y \\
\frac{dy}{dt}|_{DEV} = x
\end{cases}
\label{eq:s6}
\end{equation}

The initial conditions were sampled using a Latin-hypercube strategy, generating initial conditions for the two variables, with bounds set between [0.1, 1] for both variables.

\section{Determination of the Regularisation Parameter}
The overall regularisation term in the loss function \eqref{eq:5bis} is 
\begin{equation}
\lambda_{1,\xi} \sum_{i=0}^{N_L} \sum_{j=0}^{N_S} \left| \xi_{ij} \right| + \lambda_{1,\delta} \sum_{j=0}^{N_S} \delta_j.
\end{equation}
To ensure that both components of the regularisation term contribute on the same order of magnitude, we impose the condition
\begin{equation}
\lambda_{1,\xi} \sum_{i=0}^{N_L} \sum_{j=0}^{N_S} \left| \xi_{ij} \right| = \lambda_{1,\delta} \sum_{j=0}^{N_S} \delta_j.
\end{equation}
In the worst-case scenario, where all model parameters reach their upper or lower bounds and all terms are active, we assume
\begin{equation}
\sum_{i=0}^{N_L} \sum_{j=0}^{N_S} \left| \xi_{ij} \right| = \max(|ub|, |lb|) \cdot N_L \cdot N_S.
\end{equation}
Substituting this into the equality condition yields
\begin{equation}
\lambda_{1,\xi} \cdot \max(|ub|, |lb|) \cdot N_L \cdot N_S = \lambda_{1,\delta} \cdot N_S,
\end{equation}
which simplifies to
\begin{equation}
\lambda_{1,\delta} = \lambda_{1,\xi} \cdot \max(|ub|, |lb|) \cdot N_L.
\end{equation}
This leads to the expression provided in Equation \eqref{eq:7}.

\section{Experimental demonstration}
\subsection{Effect of the experimental noise}
\begin{figure}[H]
    \centering
    \includegraphics[width=0.9\textwidth]{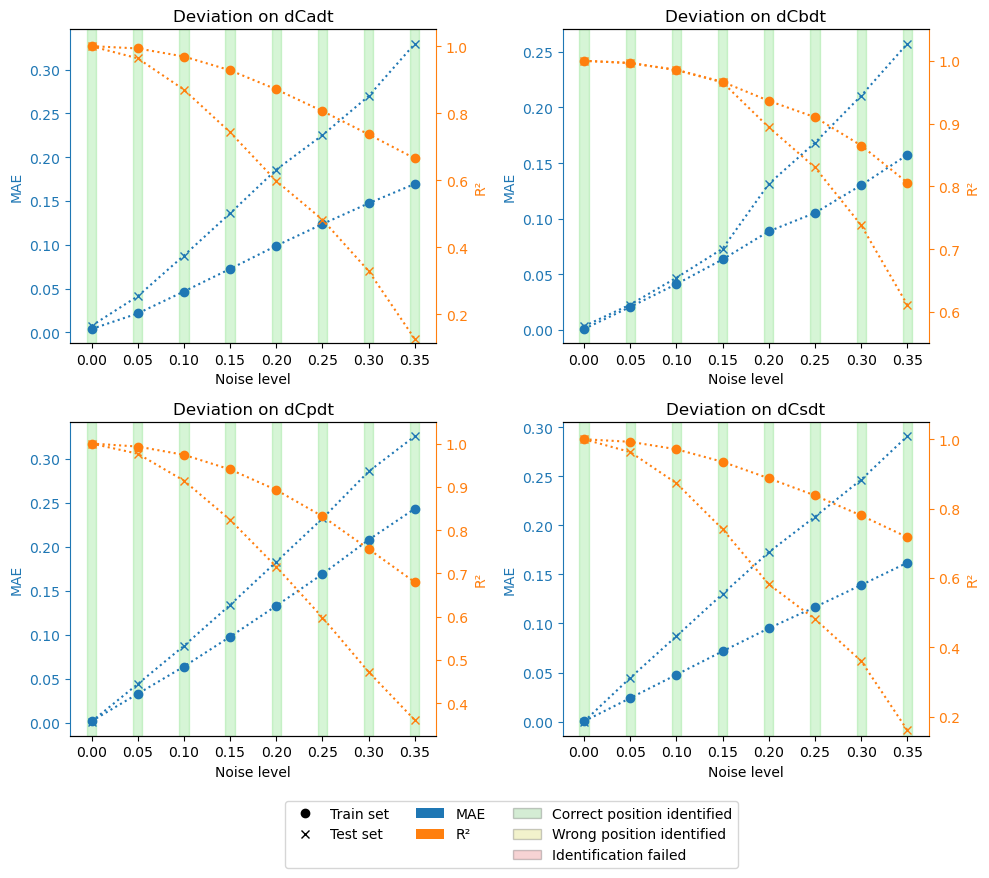}
    \caption{Effect of the experimental noise on the HM identification and training for the Meerwein arylation case, detailed with the location of epistemic uncertainty positions with an experimental noise up to 35\%. The model can identify the position of the epistemic uncertainty correctly while returning accurate models. The models worsen their predictions by increasing the experimental noise because of the presence of noise that dilutes the information to capture.}
    \label{fig:1s}
\end{figure}

\begin{figure}[H]
    \centering
    \includegraphics[width=0.9\textwidth]{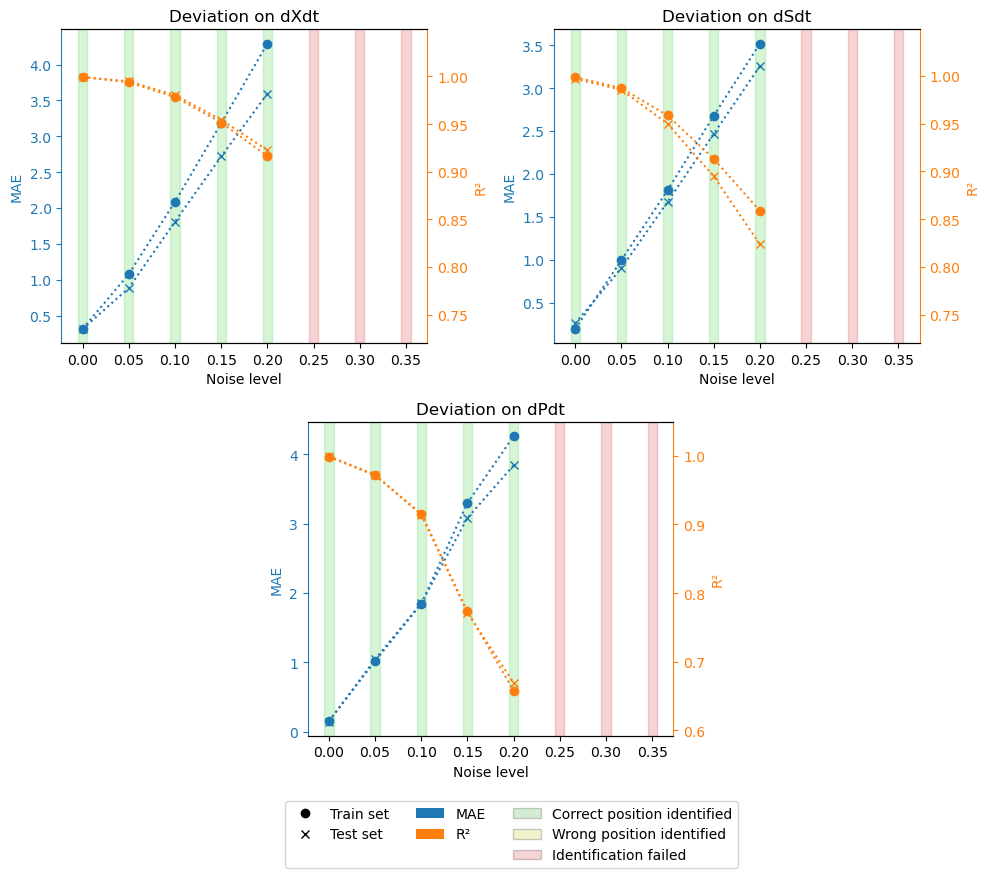}
    \caption{Effect of the experimental noise on the HM identification and training for the continuous fermentation case, detailed with the location of epistemic uncertainty positions with an experimental noise up to 35\%. The model can identify the position of the epistemic uncertainty accurately while returning accurate models only for noise up to 20\%. The models worsen their predictions by increasing the experimental noise because of the presence of noise that dilutes the information to capture. On the other hand, with noise above 20\%, the identification failed.}
    \label{fig:2s}
\end{figure}

\begin{figure}[H]
    \centering
    \includegraphics[width=0.9\textwidth]{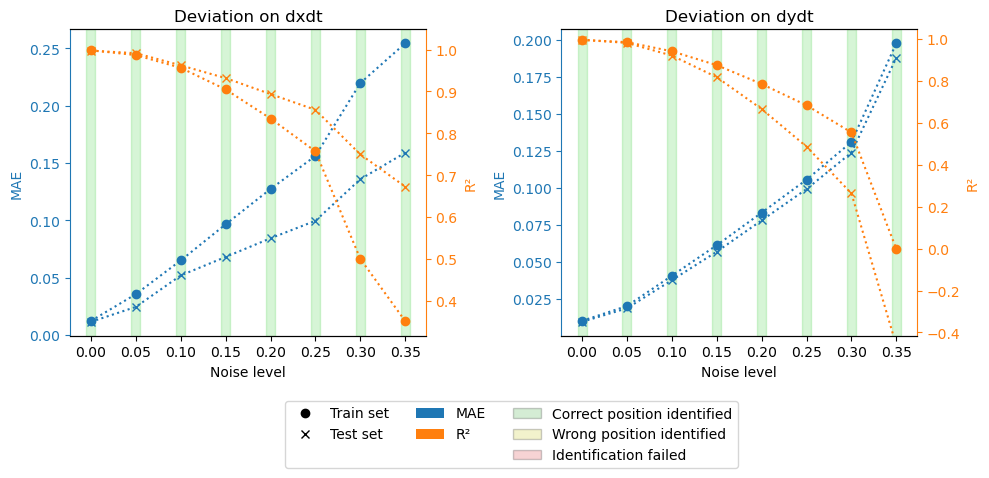}
    \caption{Effect of the experimental noise on the HM identification and training for the Lotka-Volterra validation case, detailed with the location of epistemic uncertainty positions with an experimental noise up to 35\%. The model can identify the position of the epistemic uncertainty correctly while returning accurate models. The models worsen their predictions by increasing the experimental noise because of the presence of noise that dilutes the information to capture.}
    \label{fig:3s}
\end{figure}

\subsection{Effect of the regularisation factor}
\begin{figure}[H]
    \centering
    \includegraphics[width=0.9\textwidth]{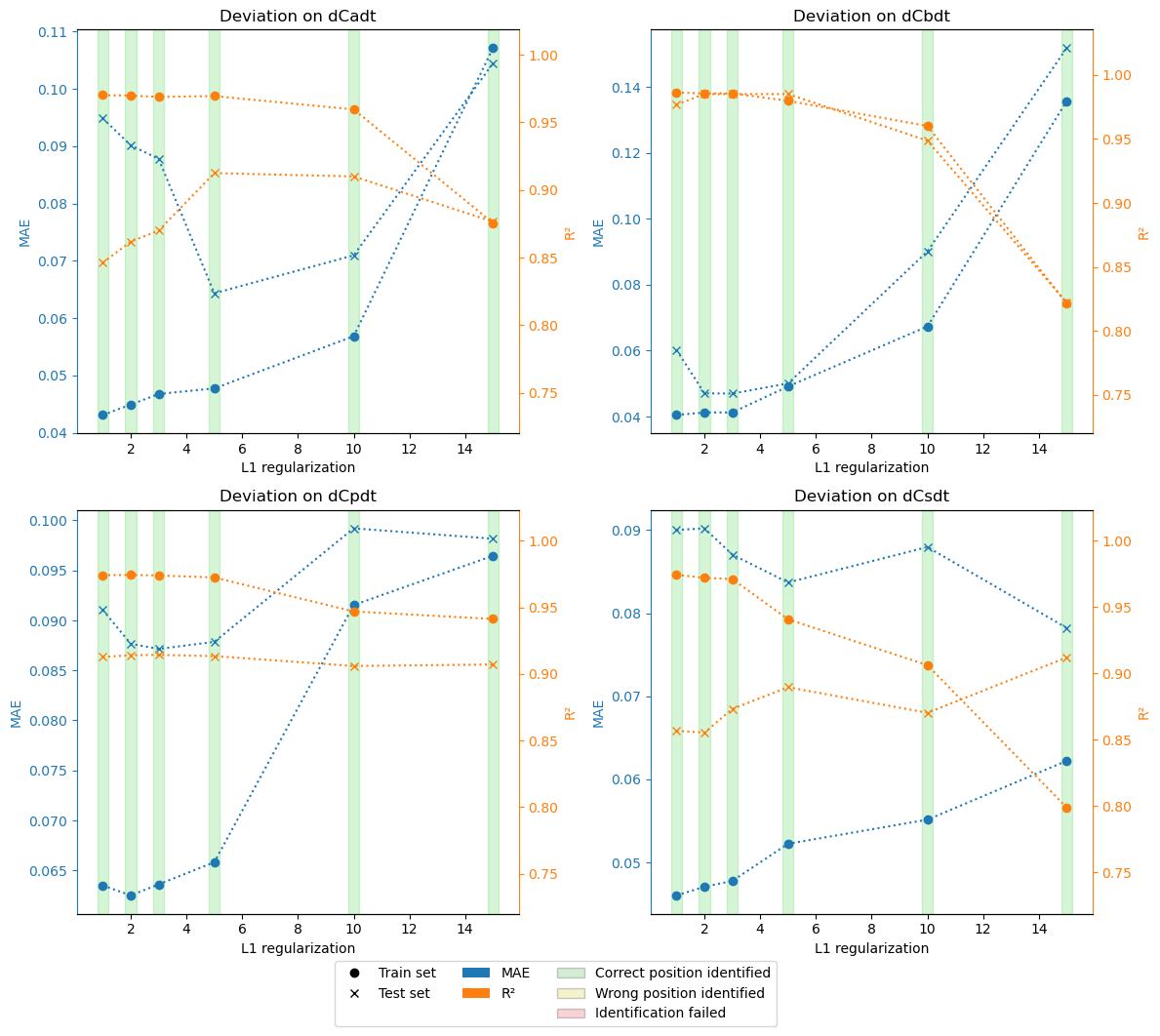}
    \caption{Effect of the L1 regularisation factor on the Meerwein arylation case. The model performance on training and testing sets worsens for $\lambda_1$ values above 3 in all the applied deviations}
    \label{fig:4s}
\end{figure}

\begin{figure}[H]
    \centering
    \includegraphics[width=0.9\textwidth]{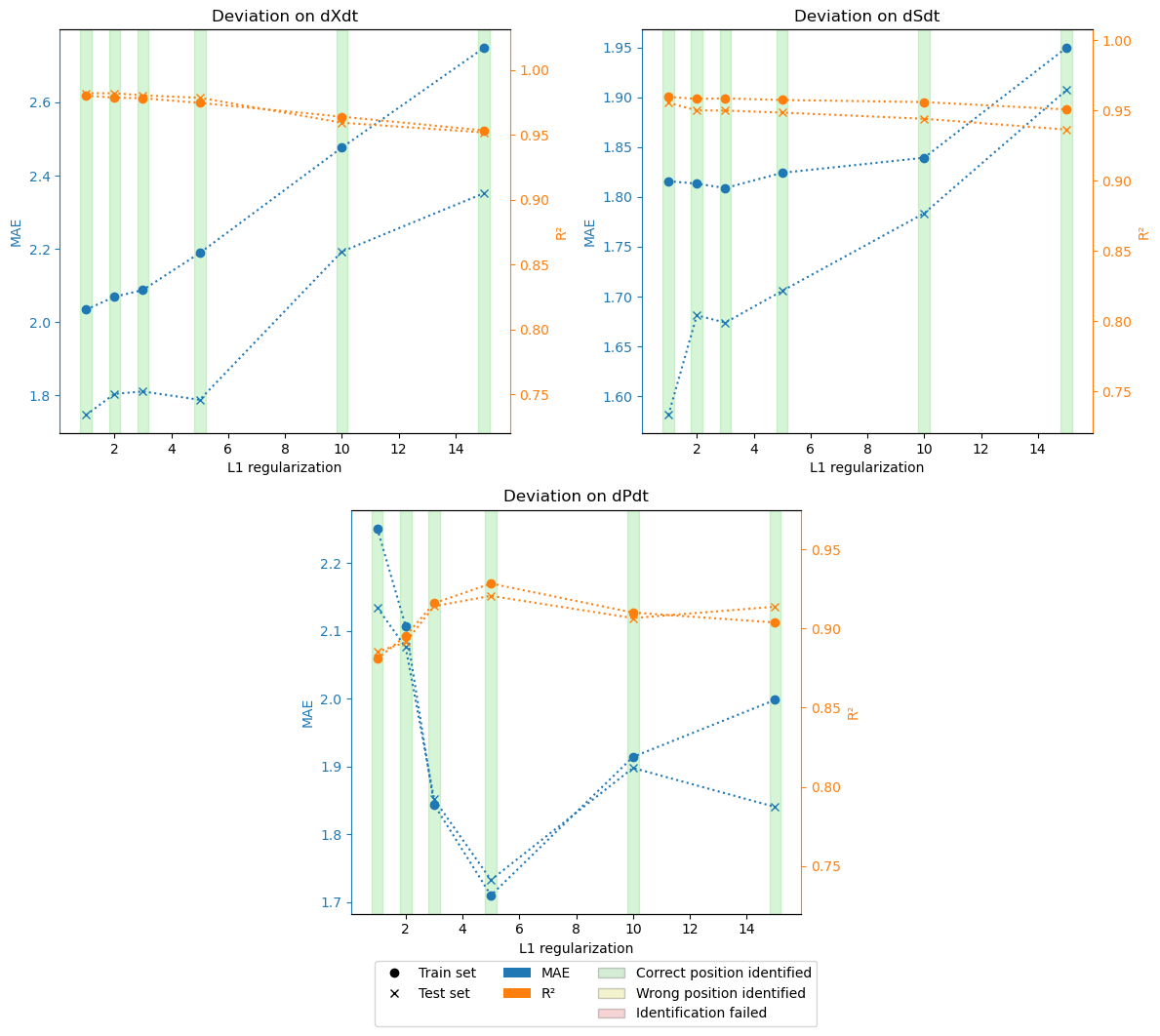}
    \caption{Effect of the L1 regularisation factor on the continuous fermentation case. The performance of the model on training and testing sets worsens for $\lambda_1$ values above 3 in all the applied deviations.}
    \label{fig:5s}
\end{figure}

\begin{figure}[H]
    \centering
    \includegraphics[width=0.9\textwidth]{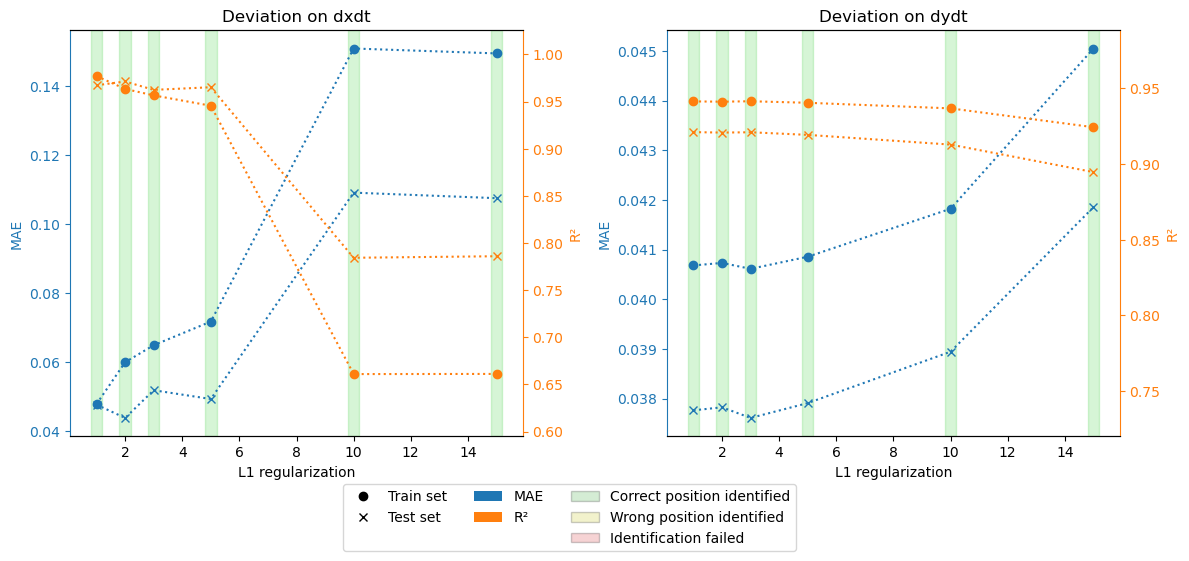}
    \caption{Effect of the L1 regularisation factor on the Lotka-Volterra oscillator case. The performance of the model on training and testing sets worsens for $\lambda_1$ values above 3 in all the applied deviations.}
    \label{fig:6s}
\end{figure}

\subsection{Effect of the time-samples amount}
\begin{figure}[H]
    \centering
    \includegraphics[width=0.9\textwidth]{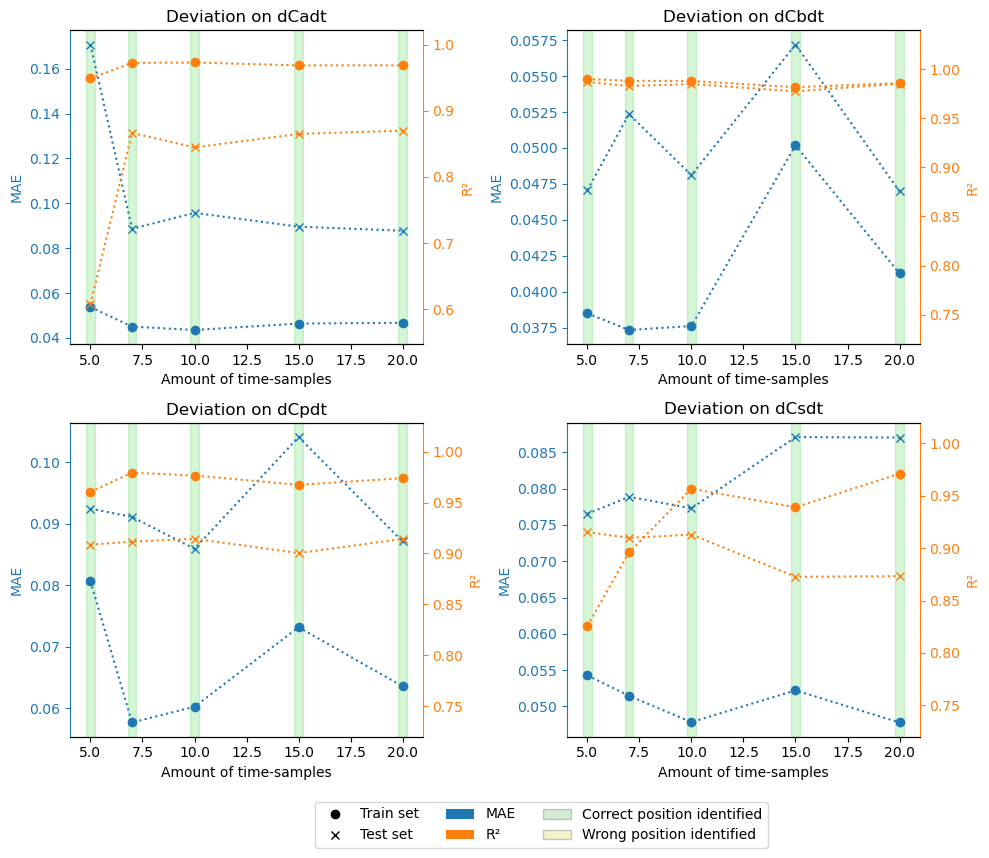}
    \caption{Effect of the number of time samples per experiment for the runs included in the training set for the Meerwein arylation case, detailed with the location of epistemic uncertainty positions.}
    \label{fig:7s}
\end{figure}

\begin{figure}[H]
    \centering
    \includegraphics[width=0.9\textwidth]{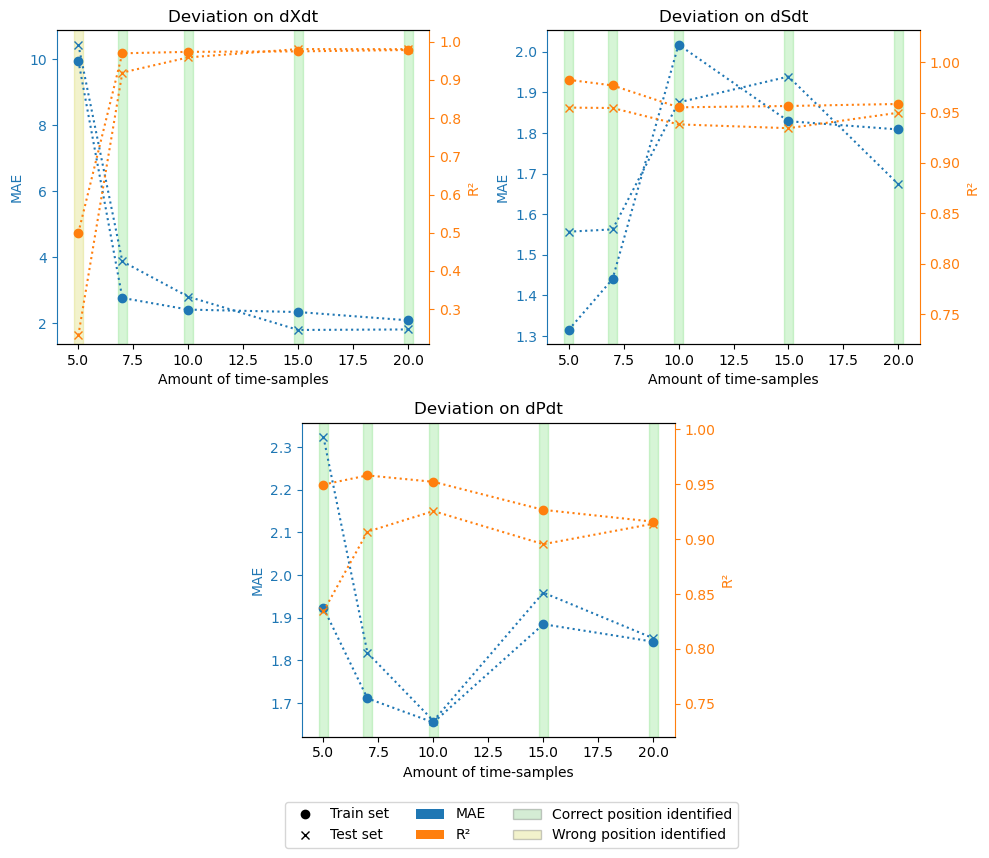}
    \caption{Effect of the number of time samples per experiment for the runs included in the training set for the continuous fermentation case, detailed with the location of epistemic uncertainty positions.}
    \label{fig:8s}
\end{figure}

\begin{figure}[H]
    \centering
    \includegraphics[width=0.9\textwidth]{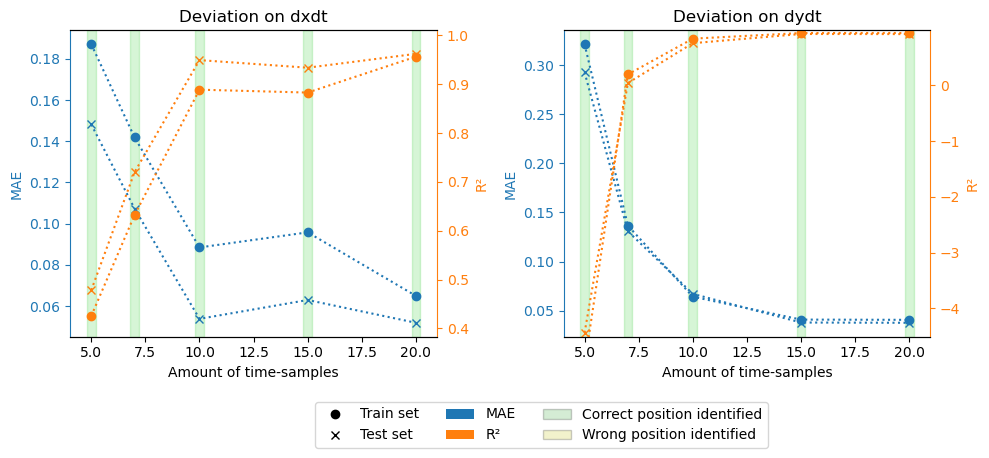}
    \caption{Effect of the number of time samples per experiment for the runs included in the training set for the Lotka-Volterra validation case detailed with the location of epistemic uncertainty positions.}
    \label{fig:9s}
\end{figure}

\subsection{Effect of the deviation on distribution}
\begin{figure}[H]
    \centering
    \includegraphics[width=0.9\textwidth]{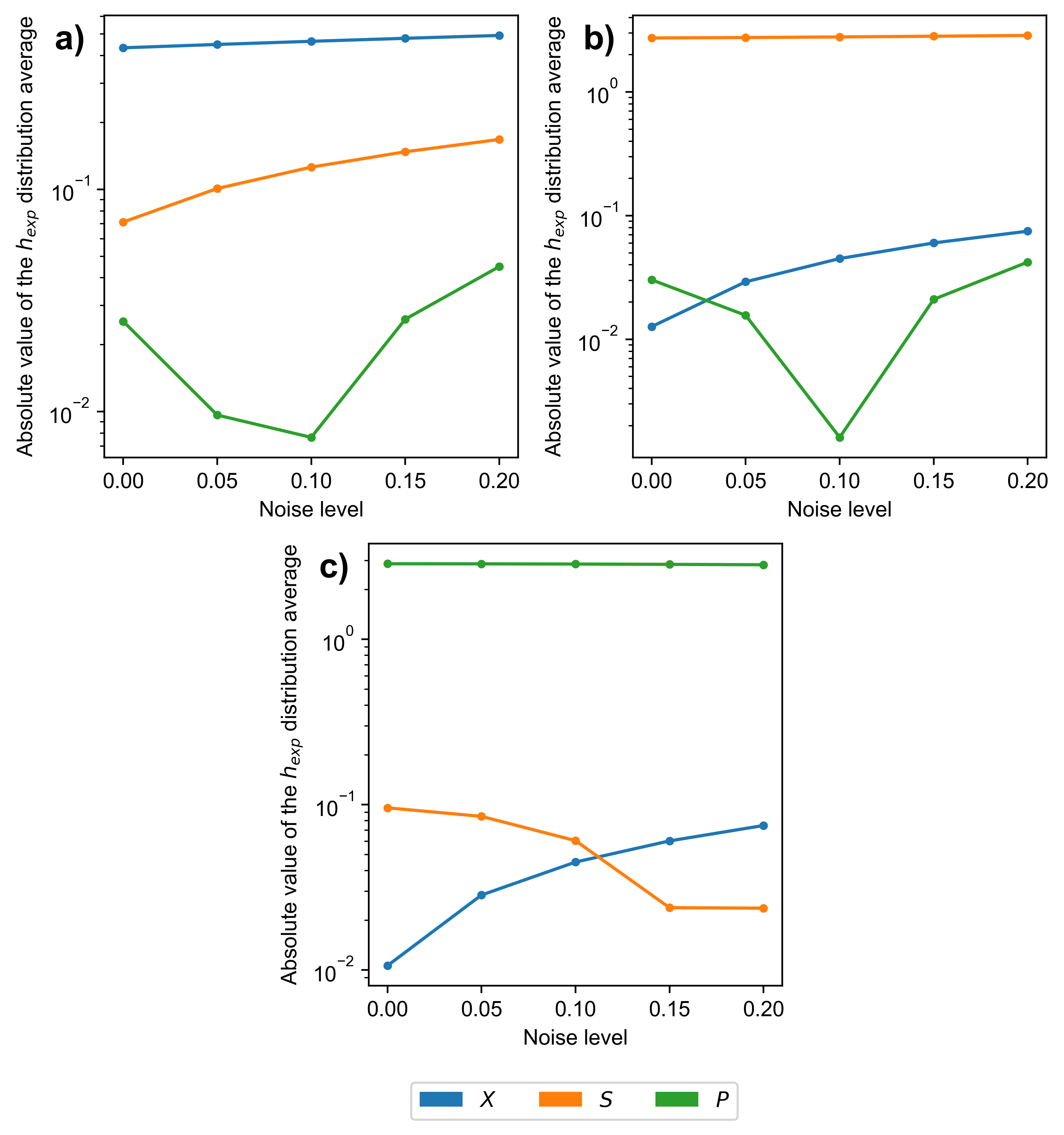}
    \caption{Values of the distribution average for the system deviation in the $h$ matrix at various noise levels for the continuous fermentation case. a) Epistemic uncertainty on $X$ dynamic, b) Epistemic uncertainty on $S$ dynamic, c) Epistemic uncertainty on $P$ dynamic}
    \label{fig:10s}
\end{figure}

\begin{figure}[H]
    \centering
    \includegraphics[width=0.9\textwidth]{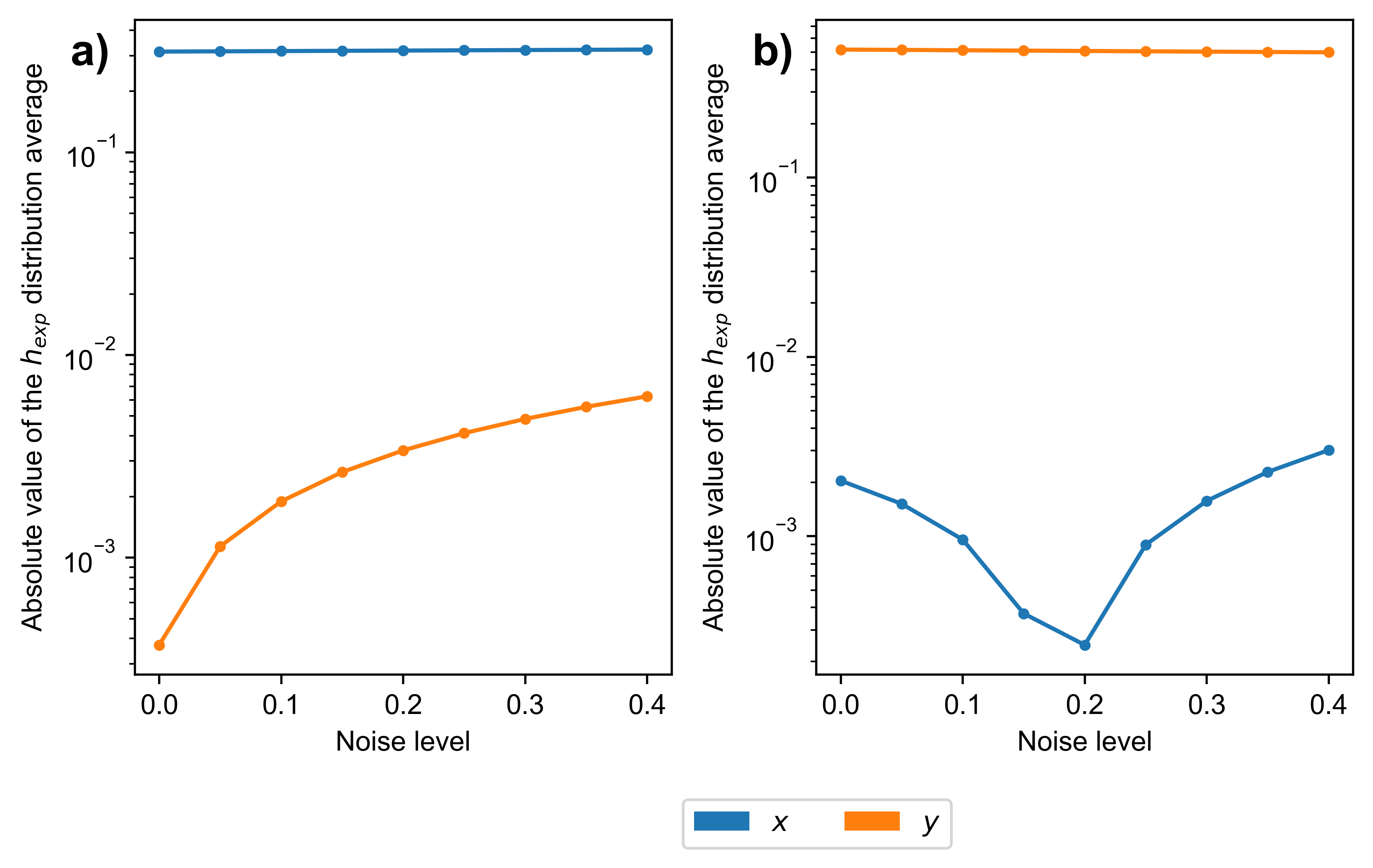}
    \caption{Values of the distribution average for the system deviation in the $h$ matrix at various noise levels for the Lotka-Volterra oscillator case. a) Epistemic uncertainty on $x$ dynamic, b) Epistemic uncertainty on $y$ dynamic}
    \label{fig:11s}
\end{figure}

\end{document}